%% file: Transfer-Learning-for-Mean-Function-Estimation.tex
\documentclass[aos, addressatend]{imsart}

\RequirePackage{amsthm,amsmath,amsfonts,amssymb}
\RequirePackage[numbers]{natbib}
\RequirePackage[colorlinks,citecolor=blue,urlcolor=blue]{hyperref}
\RequirePackage{graphicx, xcolor}
\RequirePackage[bb=dsserif]{mathalpha}
\RequirePackage{mathtools}

\startlocaldefs
\theoremstyle{plain}
\newtheorem{theorem}{Theorem}[section]

\newtheorem{proposition}[theorem]{Proposition}
\theoremstyle{remark}
\newtheorem{assumption}{Assumption}

\endlocaldefs

\input{macro.sty}
\usepackage{subcaption}
\usepackage{algorithm,algpseudocode,tabularx}

\makeatletter
\newcommand{\multiline}[1]{%
  \begin{tabularx}{\dimexpr\linewidth-\ALG@thistlm}[t]{@{}X@{}}
    #1
  \end{tabularx}
}
\makeatother

\begin{document}

\input{content/front.tex}
\input{content/introduction.tex}
\input{content/common-design.tex}
\input{content/independent-design.tex}
\input{content/numerical-experiments.tex}
\input{content/discussions.tex}
\input{content/proofs.tex}
\input{content/epilogue.tex}

\bibliographystyle{imsart-nameyear} 
\bibliography{reference}   

\end{document}

%% file: content/front.tex

\begin{frontmatter}
    \title{Transfer Learning for Functional Mean Estimation: Phase Transition and Adaptive Algorithms}
    \runtitle{Transfer Learning for Functional Mean Estimation}
    \begin{aug}
        \author{\fnms{T. Tony}~\snm{Cai}\ead[label=e1]{tcai@wharton.upenn.edu}
        \ead[label=u1,url,text=http://www-stat.wharton.upenn.edu/\textasciitilde tcai/]{http://www-stat.wharton.upenn.edu/~tcai/}},
        \author{\fnms{Dongwoo}~\snm{Kim}\(^\dag\)\ead[label=e2]{dongwooo@wharton.upenn.edu}},
        \and
        \author{\fnms{Hongming}~\snm{Pu}\(^\dag\)\ead[label=e3]{hpu@wharton.upenn.edu}}
        \runauthor{T. T. Cai, D. Kim, and H. Pu}
        \address{Department of Statistics and Data Science \\
            The Wharton School, University of Pennsylvania \\
            265 South 37th Street \\
            Philadelphia, Pennsylvania 19104 \\
            USA \\
            E-mail: \printead{e1} \\ 
            \phantom{E-mail:} \printead{e2} \\ 
            \phantom{E-mail:} \printead{e3} \\ 
            URL: \printead{u1} \\ 
        }
    \end{aug}

    \begin{abstract}
    This paper studies transfer learning for estimating the mean of random functions based on discretely sampled data, where, in addition to observations from the target distribution, auxiliary samples from similar but distinct source distributions are available. The paper considers both common and independent designs and establishes the minimax rates of convergence for both designs. The results reveal an interesting phase transition phenomenon under the two designs and demonstrate the benefits of utilizing the source samples in the low sampling frequency regime.

    For practical applications, this paper proposes novel data-driven adaptive algorithms that attain the optimal rates of convergence within a logarithmic factor simultaneously over a large collection of parameter spaces. The theoretical findings are complemented by a simulation study that further supports the effectiveness of the proposed algorithms.
    \end{abstract}



    \def\thefootnote{\(\dag\)}\footnotetext{These authors contributed equally to this work.}\def\thefootnote{\arabic{footnote}}
\end{frontmatter}


%% file: content/introduction.tex

\section{Introduction}
  
Functional data is commonly observed in a wide range of fields, and there has been significant research on functional data analysis (FDA) with applications in diverse areas such as biomedical studies, neuroscience, linguistics, psychology, demography, economics, and engineering. For a comprehensive review of the FDA and its applications, we recommend referring to \citet{ramsayAppliedFunctionalData2002} and \citet{wangFunctionalDataAnalysis2016}.

Estimating the mean function is a fundamental problem in the FDA that has garnered considerable attention in the existing literature. Notable references include \citet{riceEstimatingMeanCovariance1991, pageNormalizingTemporalPatterns2006, parkClassificationGeneFunctions2008, jiangSmoothingDynamicPositron2009}, and \citet{caiOptimalEstimationMean2011}. This problem arises naturally in a broad range of practical applications, such as analyzing a diffusion tensor imaging (DTI) dataset to investigate multiple sclerosis (MS) patients \citep{staicuModelingFunctionalData2012, degrasSimultaneousConfidenceBands2017, pomannTwoSampleDistributionFree2016}, analyzing a longitudinal CD4 cell count dataset to investigate acquired immunodeficiency syndrome (AIDS) \citep{zhouEfficientEstimationNonparametric2018}, and examining a longitudinal dataset of trajectories of patient-reported symptom severity of Parkinson's disease \citep{weaverFunctionalDataAnalysis2021}.

On the other hand, transfer learning has experienced growing popularity as a machine learning technique aimed at improving performance in a target domain by utilizing information from different yet related source domains. This approach proves particularly valuable in situations where acquiring target observations is infrequent or costly, but data from analogous studies are available. Transfer learning has been applied in various machine learning applications, including computer vision \citep{tzengAdversarialDiscriminativeDomain2017,gongGeodesicFlowKernel2012}, speech recognition \citep{huangCrosslanguageKnowledgeTransfer2013}, and genre classification \citep{choiTransferLearningMusic2017}. We refer readers to \citet{panSurveyTransferLearning2010} and \citet{weissSurveyTransferLearning2016} for a more comprehensive discussion on the practical applications of transfer learning. Recently, transfer learning has also attracted increasing attention in statistics and has been successful in various statistical learning problems such as nonparametric classification \citep{caiTransferLearningNonparametric2021, reeveAdaptiveTransferLearning2021}, high-dimensional linear regression \citep{liTransferLearningHighDimensional2022}, large-scale Gaussian graphical model \citep{liTransferLearningLargeScale2022}, high-dimensional generalized linear models \citep{liEstimationInferenceHighDimensional2023,tianTransferLearningHighDimensional2022}, nonparametric regression \citep{caiTransferLearningNonparametric2022}, and contextual multi-armed bandits \citep{caiTransferLearningContextual2022}.

Transfer learning can also be effective in the context of functional mean estimation when subjects are divided into groups, such as medical measurement datasets grouped by disease status \citep[Chapter 5]{ramsayAppliedFunctionalData2002}, gene expression datasets grouped by particular disease \citep{lengClassificationUsingFunctional2006, songOptimalClassificationTimecourse2008, parkClassificationGeneFunctions2008, wuFunctionalEmbeddingClassification2010}, or spatial dispersion datasets of marine mammals grouped by species \citep{manteFunctionalDataanalyticApproach2005}. The application of transfer learning can improve the performance of functional mean estimation in the target group by utilizing knowledge from other groups.

An interesting and motivating example of functional mean estimation with transfer learning is found in the trend analysis of coronavirus disease-19 (COVID-19). The objective is to gain insights into the epidemic's progression and identify overall trends in the United States. Similar to the work by \citet{kalogridisRobustOptimalEstimation2022}, one can look at the number of new daily COVID-19 cases across the 50 U.S. states. Each state's curve may be regarded as a functional realization and the problem reduces to estimate their functional mean. However, COVID-19 has had a global impact, affecting populations worldwide \citep{kozloffCOVID19GlobalPandemic2020}. Structural similarities in trends across different countries have been observed, as demonstrated by \citet{jamesTrendsCOVID19Prevalence2021}. This suggests that transfer learning can improve the estimation accuracy of the COVID-19 trend in the United States by utilizing the COVID-19 curves from other American countries or other continents, including Africa, Asia, Europe, and Oceania.

\subsection{Problem formulation}

The problem of functional mean estimation in a conventional setting can be described as follows. Let \(X^{[t]}:[0,1] \to \Rbb\) be a target random function and each target subject \(i \in \{1,\ldots,n^{[t]}\}\) has an an independent copy \(X^{[t]}_i\) (called curve) of \(X^{[t]}\). We have noisy observations of these curves at discrete locations:
\[ Y_{ij}^{[t]}=X_{i}^{[t]}(T_{ij}^{[t]})+\epsilon_{ij}^{[t]}, \quad(j=1,\ldots, m_i^{[t]} ~\text{and}~ i=1,\ldots, n^{[t]}), \numberthis\label{eq:modeling for target sample} \]
where \(T^{[t]}_{ij}\) are target design points, and \(\epsilon^{[t]}_{ij}\) are an independent random noises. Based on the target sample \(\Dcal^{[t]}:= \{(T_{ij}^{[t]}, Y_{ij}^{[t]}): j=1,\ldots, m_i^{[t]},~ i=1,\ldots, n^{[t]}\}\), our primary objective is to estimate the target mean function \(f^{[t]}(\cdot):=\Ebb(X^{[t]}(\cdot))\). Such problems naturally arise in a wide range of applications and are typical in FDA. Extensive examples can be found in the works of \citet{ramsayAppliedFunctionalData2002}.

In the transfer learning setup, we also have \(K\) source samples \(\Dcal^{[s,1]}, \ldots, \Dcal^{[s,K]}\) in addition to the target sample \(\Dcal^{[t]}\). These source samples are generated similarly. That is, for each source index \(k \in \{1,\ldots,K\}\), there is a random function \(X^{[s,k]}:[0,1]\to\Rbb\) with its source mean function \(f^{[s,k]}(\cdot):=\Ebb(X^{[s,k]}(\cdot))\) and we observe 
\[ Y_{ij}^{[s,k]}=X_{i}^{[s,k]}(T_{ij}^{[s,k]})+\epsilon_{ij}^{[s,k]}, ~~(j=1,\ldots, m_{i}^{[s,k]},~ i=1,\ldots, n^{[s,k]}, ~\text{and}~ k=1, \ldots, K), \numberthis\label{eq:modeling for source samples} \]
where the curves \(X_i^{[s,k]}\) are an independent copies of \(X^{[s,k]}\), \(T_{ij}^{[s,k]}\) are source design points, and \(\epsilon_{ij}^{[s,k]}\) are an independent random noises. The estimator for the target mean function \(f^{[t]}\) now can utilize the source samples \(\Dcal^{[s,1]}, \ldots, \Dcal^{[s,K]}\) as well as the target sample \(\Dcal^{[t]}\).

We model the relationship between \(\Dcal^{[s,1]}, \ldots, \Dcal^{[s,K]}\) and \(\Dcal^{[t]}\) through their corresponding mean functions \(f^{[s,1]},\ldots,f^{[s,K]}\) and \(f^{[t]}\). Let us denote by \(\Hcal_{\alpha}(L, M)\) the class of bounded functions with H\"{o}lder smoothness \(\alpha > 0\). To be more specific, \(f \in \Hcal_{\alpha}(L, M)\) if and only if \(f:[0,1]\rightarrow \Rbb\) is bounded as \(\norm{f}_{\Lcal^\infty} \leq M\) and \(\alpha^*\)-times continuously differentiable such that
\[ \abs{f^{(\alpha^*)}(t_1)-f^{(\alpha^*)}(t_2)} \leq L \abs{t_1-t_2}^{\alpha-\alpha^*}, \qquad\text{for any}~ t_1,t_2\in [0,1], \] 
where \(\alpha^* := \omega(\alpha)\) denotes the largest integer strictly smaller than \(\alpha\). We assume 
\[\begin{aligned}
&f^{[t]}, f^{[s,k]} \in \Hcal_{\alpha_m}(L_m, M_m), \\
&\delta^{[s,k]} := f^{[t]} - f^{[s,k]} \in \Hcal_{\alpha_\delta}(L_\delta, M_\delta),
\end{aligned} \qquad(k=1,\ldots,K) \numberthis\label{eq:modeling assumption for mean and difference functions}\]
with two smoothness parameters \(\alpha_m, \alpha_\delta > 0\) and constants \(L_m, M_m, L_\delta, M_\delta > 0\). Our model captures a broader range of flexibility and provides more informative insights into the problem. We have never made any specific assumptions or constraints regarding the relationships between smoothness parameters, \(\alpha_m\) and \(\alpha_\delta\). This allows us to adapt to a wider variety of scenarios and better address the complexities of the problem.

In line with the conventional framework for functional mean estimation \citep{caiOptimalEstimationMean2011}, we consider two distinct sampling designs. The first design is characterized by common design points, where observations are gathered at identical locations across curves. In this setting, we have \(T_{ij}^{[t]} = T_{j}^{[t]}\) and \(T_{ij}^{[s,k]} = T_{j}^{[s]}\). The second design employs an independent approach, where \(T_{ij}^{[t]}\) and \(T_{ij}^{[s,k]}\) are independently sampled from a distribution defined over \([0,1]\). For a more comprehensive formalization of these designs, Section \ref{sec:optimality under a common design} is dedicated to exploring a common design setting, while Section \ref{sec:optimality under an independent design} delves into an independent design case.

\subsection{Our contribution}

To streamline our presentation, we make an assumption that the source samples share the same design type as the target sample. Additionally, we assume that both the target and source samples exhibit identical characteristics, including the number of subjects and design points per subject. In other words,
\[\begin{aligned}
    &n^{[t]} = n_t ~\text{and}~ m_i^{[t]}=m_t &&\text{for any}~ i=1,\ldots,n^{[t]}, \\
    &n^{[s,k]} = n_s ~\text{and}~ m_i^{[s,k]}=m_s &&\text{for any}~ i=1,\ldots,n^{[s, k]} ~\text{and}~ k=1,\ldots,K.
\end{aligned} \numberthis\label{eq:modeling assumption for the number of subjects and design points}\]
We propose novel algorithms and develop an optimality theory for the estimation of the target mean function \(f^{[t]}\) under both common and independent designs. Under a common design, the optimal minimax rate of convergence up to logarithmic factors is shown to be
\[ \left[m_t^{-2\alpha_m} + \frac{1}{n_t}\right] \wedge \left[m_t^{-2\alpha_\delta} + \frac{1}{n_t} + m_s^{-2\alpha_m} + \frac{1}{K n_s}\right]. \]
On the other hand, under an independent design, the rate is given by
\[\left[(m_t n_t)^{-2\alpha_m/(2\alpha_m+1)} + \frac{1}{n_t}\right] \binmin \left[(m_t n_t)^{-2\alpha_\delta/(2\alpha_\delta+1)} + \frac{1}{n_t} + (K m_s n_s)^{-2\alpha_m/(2\alpha_m+1)} + \frac{1}{K n_s} \right].\]

The examination of minimax convergence rates uncovers several critical aspects of phase transition that significantly impact the effectiveness of transfer learning. First of all, for any model in which \(\alpha_\delta \leq \alpha_m\) is satisfied, the utilization of source samples and transfer learning becomes futile. Under such a model, the convergence rates are identical to those achieved when \(n_s = 0\). Essentially, any model with \(\alpha_\delta \leq \alpha_m\) cannot outperform the conventional learning setup. The effectiveness of transfer learning is justifiable only for models such that \(\alpha_\delta\) is strictly greater than \(\alpha_m\). This observation gives rise to the phase transition phenomenon, highlighting the importance of modeling assumption \(\alpha_\delta > \alpha_m\) for meaningful discussions regarding the effectiveness of transfer learning.

Under the extra modeling assumption \(\alpha_\delta > \alpha_m\), we can find crucial differences between two regimes: high and low sampling frequencies. In the regime of low sampling frequency where \(m_t \ll n_t^{1/2\alpha_m}\), transfer learning can lead to faster convergence rates compared to the case of utilizing the target sample only. However, in the high sampling frequency regime where the target sampling frequency \(m_t\) grows at least as rapidly as the rate \(n_t^{1/2\alpha_m}\), transfer learning cannot improve the performance beyond the optimal parametric rate of \(n_t^{-1}\), resulting in a phase transition at \(m_t \asymp n_t^{1/2\alpha_m}\). The remarkable finding is that this transition boundary remains invariant, regardless of whether the design is common or independent, which underscores the structural equivalence between the two designs.
 
We turn our attention to the low sampling frequency regime where \(m_t \ll n_t^{1/2\alpha_m}\) holds. In this regime, we shall explore distinct conditions under which transfer learning can enhance performance under common and independent designs, respectively. Under a common design, transfer learning proves effective when \(m_s^{-{2\alpha_m}}+ (Kn_s)^{-1}\ll m_t^{-{2\alpha_m}}\), whereas it is effective if \((K m_s n_s)^{-2\alpha_m/(2\alpha_m+1)}+(Kn_s)^{-1} \ll (m_t n_t)^{-2\alpha_m/(2\alpha_m+1)}\) under an independent design. To a certain extent, both designs require a sufficient number of total subjects in the source samples \((K n_s \gg m_t^{2\alpha_m} ~\text{and}~ K n_s \gg (m_t n_t)^{2\alpha_m/(2\alpha_m+1)}, ~\text{respectively})\), but neither implies the other. The second condition, however, exhibits a significant disparity between the two designs. In a common design, a higher sampling frequency in the source samples (\(m_s \gg m_t\)) is necessary, while in an independent design, a larger total number of observations in the source samples (\(K m_s n_s \gg m_t n_t\)) is required. Finally, another noteworthy observation is that the size of the source group \(K\) plays a more influential role in an independent design than in a common design. This is mainly because increasing \(K\) leads to a larger total number of observations in the source samples, but not to a higher sampling frequency in the source samples.

A comparison of the minimax risks between common and independent designs reveals that the rate of convergence for an independent design is consistently faster or equal to that of a common design. Consequently, when confronted with the choice between these two designs, assuming all other factors remain constant, selecting an independent design should be always preferred. This observation aligns with our intuition that estimating the functional mean will be more accurate when approached holistically rather than atomistically. In this regard, an independent design possesses inherent advantages as it explores a wider range of design points than a common design does, even though the accuracy of estimation for each design point may be relatively lower.

Although our primary focus is on transfer learning, these results also provide new insights into functional mean estimation in the conventional setup. By setting both the source sample size \(n_s\) and sampling frequency \(m_s\) to zero, we obtain rates for functional mean estimation without transfer learning under both common and independent designs. Our rates of convergence, \(m_t^{-2\alpha_m}+n_t^{-1}\) and \((m_t n_t)^{-2\alpha_m/(2\alpha_m+1)}+n_t^{-1}\), serve as a generalization of the rates derived by \citet{caiOptimalEstimationMean2011}. This generalization is particularly significant when the target mean function is rough and even not differentiable, as the current framework still guarantees consistent estimation. Such flexibility is a key advantage of our model.

In this paper, we introduce two novel algorithms, \(\Acal_{\mathrm{CL}}\) and \(\Acal_{\mathrm{TL}}\), where the latter utilizes the source samples to transfer knowledge but the former does not. These algorithms can be applied under both designs and achieve the optimal rate with carefully chosen parameters. However, these optimal parameters differ between the two designs. Under an independent design, the optimal values depend on the unknown modeling parameters \(\alpha_m\) and \(\alpha_\delta\) while in a common design, they do not. Therefore, for constructing an adaptive and optimal procedure under a common design, it suffices to incorporate the two algorithms \(\Acal_{\mathrm{CL}}\) and \(\Acal_{\mathrm{TL}}\), while under an independent design, one additionally needs to tune parameters adaptively. This highlights the difference between the two designs in constructing adaptive estimators. We propose adaptive procedures \(\Acal_{\mathrm{ALC}}\) and \(\Acal_{\mathrm{ALI}}\) for common and independent designs, respectively, and demonstrate their optimality.

\subsection{Organization and notation}

The rest of the paper is organized as follows. We finish this section with basic notation and then study transfer learning for functional mean estimation under a common design in Section \ref{sec:optimality under a common design}. The optimal rate of convergence is established and an adaptive algorithm is proposed. We next shift our focus to an independent design in Section \ref{sec:optimality under an independent design}. In Section \ref{sec:numerical experiments}, we conduct a simulation study to investigate the numerical performance of the proposed adaptive algorithms for both common and independent designs. The numerical results agree with the theoretical findings presented earlier. Section \ref{sec:discussion} discusses further research directions. The proofs of the main results are provided in Section \ref{sec:proof}, while additional proofs are presented in the Supplementary Material \citep{caiSupplementTransferLearning2024}.

Throughout the paper, we consider subject size \((n_t, n_s)\), sampling frequency \((m_t, m_s)\) and the number of source samples \((K)\) as the primary asymptotic components while all other quantities are treated as constants. We adhere to the standard for big-Oh\((O)\), big-Omega\((\Omega)\) and big-Theta\((\Theta)\) notations. On occasion, for brevity, we use the symbols \(\lesssim\), \(\gtrsim\) and \(\asymp\) as replacements for big-Oh\((O)\), big-Omega\((\Omega)\) and big-Theta\((\Theta)\), respectively. Similarly, we follow the convention for the little-Oh\((o)\) notation and briefly express \(a \ll b\) or \(b \gg a\) if and only if \(a = o(b)\) is true. Although the constant term is only allowed in the standard big-Oh\((O)\), big-Omega\((\Omega)\) and big-Theta\((\Theta)\) notations, their tilde-variants, big-Oh-tilde\((\widetilde{O})\), big-Omega-tilde\((\widetilde{\Omega})\) and big-Theta-tilde\((\widetilde{\Theta})\) accept logarithmic polynomial terms of \(n_t, n_s, m_t, m_s\) and \(K\). For any function \(f:[0,1] \to \Rbb\), it is straightforward to define functional \(\Lcal^2\)-norm, \(\norm{f}_{\Lcal^2(I)}\), and \(\Lcal^\infty\)-norm, \(\norm{f}_{\Lcal^\infty(I)}\), restricted on some interval \(I \subset [0,1]\). When \(I\) stands for the whole domain \([0,1]\), we follow the convention of leaving out \(I\) in the norm notation.

%% file: content/common-design.tex

\section{Transfer learning for functional mean estimation under a common design} \label{sec:optimality under a common design}

In this section, we will explore transfer learning for functional mean estimation under a common design, where curves are observed at the same design points within the same target or source subject groups. The goal is to introduce a novel algorithm for estimating the functional mean within the transfer learning framework and derive an upper bound for its integrated mean squared error (IMSE). We then establish the minimax rate of convergence and the optimality of the proposed algorithm by obtaining a matching lower bound. Finally, we introduce an adaptive and optimal algorithm that can be utilized in practical applications without requiring knowledge of modeling parameters.

While the results presented in this section are valid for randomly selected common design points, for the sake of clarity in our exposition, we simplify by assuming that these design points are deterministic. Additionally, we adopt the convention, without loss of generality, that the common design points are arranged in ascending order. To summarize, a common design is characterized by two collections of common design points, \(\{T_j^{[t]}:j=1,\ldots,m_t\}\) and \(\{T_j^{[s]}:j=1,\ldots,m_s\}\), such that \(T_{j_1}^{[t]} \leq T_{j_2}^{[t]}\) and \(T_{j_1}^{[s]} \leq T_{j_2}^{[s]}\) for any \(j_1 \leq j_2\) as well as the following holds:
\begin{align*}
    T_{ij}^{[t]} &= T_{j}^{[t]} &&\text{for all}~ j = 1,\ldots,m_t ~\text{and}~ i=1,\ldots,n_t, \\
    T_{ij}^{[s,k]} &= T_{j}^{[s]} &&\text{for all}~ j = 1,\ldots,m_s,~ i=1,\ldots,n_s ~\text{and}~ k = 1,\ldots,K.
\end{align*}

To describe the optimality of functional mean estimation, we consider a statistical model denoted by \(\Pcal\) which is a collection of certain probability measures. The mean and difference functions of measures in \(\Pcal\) satisfy Equation \eqref{eq:modeling assumption for mean and difference functions} and the measures themselves also satisfy the uniform sub-Gaussian condition (Assumption \ref{ass:uniformly subgaussian assumption}). Finally, the collection of target and source samples, \(\{\Dcal^{[t]}, \Dcal^{[s,1]}, \ldots, \Dcal^{[s,K]}\}\), is assumed to be independent. We continue to employ the same assumptions in Section \ref{sec:optimality under an independent design} where we argue the optimality under an independent design.

\begin{assumption} \label{ass:uniformly subgaussian assumption}
The random functions and noises are uniformly sub-Gaussian variables with a positive variance proxy \(\tau^2\). As per \citet{wainwrightHighdimensionalStatisticsNonasymptotic2019}, we assume that for every \(u > 0\), the following holds: for any given \(x \in [0,1]\) and \(k=1,\ldots,K\),
\[ \left.\begin{aligned}
    &\Pbb\Bigl( \absbig{X_1^{[t]}(x) - f^{[t]}(x)} \geq u \Bigr) \binmax \Pbb\Bigl( \absbig{X_1^{[s,k]}(x) - f^{[s,k]}(x)} \geq u \Bigr) \\
    &\Pbb\Bigl( \absbig{\epsilon_{11}^{[t]}} \geq u \Bigr) \binmax \Pbb\Bigl( \absbig{\epsilon_{11}^{[s,k]}} \geq u \Bigr)
\end{aligned} \right\} \leq 2e^{-u^2/2\tau^2} \]
\end{assumption}

\subsection{Methodology and upper bound}

In this subsection, we present a novel algorithm for estimating the target mean function. Before delving into the transfer learning problem of our main interest, it is useful to revisit the conventional framework for functional mean estimation, which is a crucial component of the transfer learning algorithm. This involves estimating the target mean function \(f^{[t]}\) when no source samples are available, i.e. \(n_s = 0\). Although this problem has been extensively investigated by \citet{caiOptimalEstimationMean2011}, and the focus of this paper is on transfer learning, it is still valuable to discuss it separately.
 
\begin{theorem}[The minimax risk under conventional setup and common design] \label{thm:minimax risk under conventional setup and common design}
    Suppose no source samples are available, i.e. \(n_s = 0\) and the fixed and common design points for the target sample satisfy \(\max_{1\leq j\leq m_t+1} \bigl(T_{j}^{[t]}-T_{j-1}^{[t]}\bigr) \leq C_t/m_t\) for some constant \(C_t > 0\), where \(T_0^{[t]} = 0\) and \(T_{m_t+1}^{[t]} = 1\). Then
    \[ \inf_{\widehat{f}^{[t]}} \sup_{\Pbb\in\Pcal} \Ebb\norm{\widehat{f}^{[t]} - f^{[t]}}^2_{\Lcal^2} = \widetilde{\Theta} \left(L_m^2 m_t^{-2\alpha_m} + \frac{1}{n_t}\right), \]
    where the infimum is taken over all estimators \(\widehat{f}^{[t]} = \widehat{f}^{[t]}(\Dcal^{[t]})\) based on the target sample.
\end{theorem}

Theorem \ref{thm:minimax risk under conventional setup and common design} generalizes the minimax rate shown by \citet{caiOptimalEstimationMean2011}, \(\Theta(m_t^{-2r}+n_t^{-1})\), where both the curve \(X^{[t]}\) and mean function \(f^{[t]}\) are assumed to be \(r\) times differentiable with  \(r\in\Zbb^+\). Notably, the current framework remains consistent in estimating the mean function even when it is not smooth \((\alpha_m < 1)\), whereas the previous one does not. Moreover, our approach provides greater flexibility by not requiring the smoothness of the target curve \(X^{[t]}\). Therefore, the proposed model represents a significant improvement over the existing literature on functional mean estimation.

The primary contribution of Theorem \ref{thm:minimax risk under conventional setup and common design} lies in the estimation method, characterized by the algorithm \(\Acal\). This algorithm is equally pivotal in the transfer learning setup. At the core of the algorithm \(\Acal\), it partitions the domain \([0,1]\) into sub-intervals and performs polynomial regression on each sub-interval to estimate the mean function. Algorithm \ref{alg:randomized local polynomial regression with thresholding} outlines the step-by-step instructions for implementing this algorithm.

The purpose of introducing a randomized collection \(\Dcal_r\) in Algorithm \(\Acal\) is to ensure the independence of the observations in the collection. Unlike Algorithm \(\Acal\), when we naively take the sub-collection \(\widetilde{\Dcal}_r:= \{(T, Y) \in \Dcal: T \in I_r\}\) of all observations whose design points fall into the interval \(I_r\), they are no longer independent in general. This is because \(\widetilde{\Dcal}_r\) may contain more than one observation from the same curve depending on a common sampling scheme. This poses some technical challenges when analyzing the resulting estimator. To address this issue, we use a randomized reduction technique. The idea is to randomly select one observation within each interval \(I_r\) per curve and discard the rest. This reduction procedure ensures that the observations used for polynomial regression are independent.

We must ensure that observations are neither under-represented nor over-represented when randomly selecting them, as only one discrete observation per subject is chosen. Under an independent design, we assume that the density generating random design points is bounded both above and below by constants (See Theorem \ref{thm:minimax risk under conventional setup and independent design} and \ref{thm:upper bound under an independent design}), which solves the under- and over-representation problem. However, under a common design, we assume fixed and common design points are not too far from each other (See Theorem \ref{thm:minimax risk under conventional setup and common design} and \ref{thm:upper bound under a common design}), which only handles the under-representation issue. To address the issue of over-representation of common design, the proposed algorithm employs a (randomized) sub-collection \(\widetilde\Dcal^{(i)} \subset \Dcal^{(i)}\). The design points of \(\widetilde\Dcal^{(i)}\) are required to be spaced at a minimum separation of \(1/m\) while still encompassing those found in \(\Dcal^{(i)}\) within a maximum distance of \(1/2m\). The introduction of sub-collection \(\widetilde\Dcal^{(i)}\) is thus expected to alleviate the issue of over-representation.

\begin{algorithm}[H]
    \caption{Randomized local polynomial regression with thresholding \(\Acal(\Dcal, b, d, M)\)} \label{alg:randomized local polynomial regression with thresholding}
    \begin{algorithmic}[1] \normalsize
    \Require A collection \(\Dcal = \{(T_{ij},Y_{ij}) : j=1,\ldots,m,~ i=1,\ldots,n\}\) of observations, a bandwidth \(b \in 1/\Zbb^+\), a degree \(d \in \Zbb^+\) of polynomial, and a threshold \(M > 0\).
    \State Partition the domain \([0,1]\) into \(\inverse{b}\)-many intervals of length \(b\). We denote those intervals by \(I_r\) \((r=1,\ldots,\inverse{b})\) from left to right.
    \State Let us denote by \(\Dcal^{(i)} := \{(T_{ij},Y_{ij}) : j=1,\ldots,m\}\) the collection of observations from the same subject \(i \in \{1,\ldots,n\}\). For index \(r=1,\ldots,\inverse{b}\), we take randomized collection \(\Dcal_r := \{(\mathbfbb{T}_{i,r},\mathbfbb{Y}_{i,r}):i=1,\ldots,n\}\) of observations where \((\mathbfbb{T}_{i,r},\mathbfbb{Y}_{i,r})\) is randomly chosen from the following process:
    \If{the collection \(\Dcal\) comes from common design}
        \State \multiline{Consider any \((1/m)\)-packing and \((1/2m)\)-covering sub-collection \(\widetilde{\Dcal}^{(i)}\) of \(\Dcal^{(i)}\) where the distance between observations is computed based on design points. The random \((\mathbfbb{T}_{i,r},\mathbfbb{Y}_{i,r})\) is now uniformly chosen from \( \{(T,Y) \in \widetilde\Dcal^{(i)} : T \in I_r\}\).}
    \ElsIf{the collection \(\Dcal\) comes from independent design}
        \State The random \((\mathbfbb{T}_{i,r},\mathbfbb{Y}_{i,r})\) is uniformly chosen from \(\{(T,Y) \in \Dcal^{(i)} : T \in I_r\}\).
    \EndIf
    \State Implement the polynomial regression of degree \(d\) on each collection \(\Dcal_r\) \((r=1,\ldots,\inverse{b})\) of observations. In other words, we are enough to compute for each \(r=1,\ldots,\inverse{b}\),
    \[ (\check{a}_{r,0}, \check{a}_{r,1}, \ldots, \check{a}_{r,d}) := \argmin_{(a_0, a_1, \ldots, a_d) \in \Rbb^{d+1}} \sum_{(T,Y) \in \Dcal_r} \Biggl[Y - \sum_{s=0}^d a_s \left(\frac{T-(q-1)b}{b}\right)^{\!s} \Biggr]^2. \]
    If the solution is not available, simply take \((\check{a}_{r,0}, \check{a}_{r,1}, \ldots, \check{a}_{r,d}) = 0\).
    \State Compute the local polynomial regression estimator \(\check{f}:[0,1] \to \Rbb\) by
    \[ \check{f}(x) := \sum_{r=1}^{\inverse{b}} \sum_{s=0}^d \check{a}_{r,s} \left(\frac{x-(r-1)b}{b}\right)^{\!s} \mathbb{1}(x \in I_r) \qquad(0 \leq x \leq 1).\]
    \State Output the final estimator \(\widehat{f}:[0,1] \to \Rbb\) through thresholding:
    \[\widehat{f} := \begin{dcases}
        \check{f} &\text{if}~ \norm{\check{f}}_{\Lcal^\infty(I_r)} \leq M, \\
        0 &\text{otherwise,} \\
    \end{dcases} \quad\text{on each}~ I_r ~(r=1,\ldots,b^{-1}). \]
\end{algorithmic}
\end{algorithm}

We will now present the optimal algorithm \(\Acal_{\mathrm{CL}}(b_t, d_t, M_t)\) for estimating the mean function under the conventional setup and common design. This algorithm is essentially the same as algorithm \(\Acal\) with the following four inputs: the target sample \(\Dcal^{[t]}\), a bandwidth \(b_t \in 1/\Zbb^+\), a degree \(d_t \in \Zbb^+\) of local polynomial and a threshold \(M_t > 0\). The output \(\widehat{f}_{\mathrm{CL}}^{[t]}\) of the conventional learning algorithm \(\Acal_{\mathrm{CL}}\) indeed attains the following minimax rate of convergence in Theorem \ref{thm:minimax risk under conventional setup and common design} when input parameters are carefully selected. The detailed selection scheme will be shortly postponed to Theorem \ref{thm:upper bound under a common design}.
\[ \sup_{\Pbb\in\Pcal} \Ebb\norm{\widehat{f}_{\mathrm{CL}}^{[t]} - f^{[t]}}_{\Lcal^2}^2 \leq O\left(L_m^2 m_t^{-2\alpha_m} + \frac{\log^2 n_t}{n_t}\right). \]

Let us shift our focus back to the transfer learning setup of the primary interest. In addition to the target sample, we have access to the source samples that could be potentially beneficial for estimating the target mean function more accurately. In this setting, the transfer learning algorithm, represented by \(\Acal_{\mathrm{TL}}\), incorporates both the target and source samples and leverages their structures to achieve a more precise estimation of the target mean function. The basic idea is to decompose the target mean function \(f^{[t]}\) into the source part \(f^{[s]} := \inverse{K} \sum_{k=1}^{K} f_k^{[s]}\) and the difference part \(\delta^{[s]} := f^{[t]} - f^{[s]}\) and to estimate them separately. The step-by-step guide for implementing this algorithm is provided in Algorithm \ref{alg:transfer learning for mean function}.

\begin{algorithm}[H]
    \caption{Transfer learning for mean function \(\Acal_{\mathrm{TL}}(b_s, b_\delta, d_s, d_\delta, M_s, M_\delta)\)} \label{alg:transfer learning for mean function}
    \begin{algorithmic}[1] \normalsize
    \Require Two bandwidths \(b_s, b_\delta \in 1/\Zbb^+\), two degrees \(d_s, d_\delta \in \Zbb^+\) of polynomial, and two thresholds \(M_s, M_\delta > 0\).
    \State Execute \(\Acal(\Dcal^{[s]}, b_s, d_s, M_s)\) with the combined source sample \(\Dcal^{[s]} := \bigcup_{k=1}^K \Dcal^{[s,k]}\). The result of this algorithm is denoted as \(\widehat{f}^{[s]}\).
    \State Compute a new sample \(\Dcal^{[\delta]} := \bigl\{\bigl(T_{ij}^{[t]}, Y_{ij}^{[t]} - \widehat{f}^{[s]}(T_{ij}^{[t]})\bigr):j=1,\ldots,m_t,~ i=1,\ldots,n_t\bigr\}\).
    \State Execute \(\Acal(\Dcal^{[\delta]}, b_\delta, d_\delta, M_\delta)\). The output of this algorithm is denoted by \(\widehat{\delta}^{[s]}\).
    \State Output our final estimator \(\widehat{f}^{[t]} := \widehat{f}^{[s]} + \widehat{\delta}^{[s]}\).
\end{algorithmic}
\end{algorithm}

It is important to notice that the transfer learning algorithm \(\Acal_{\mathrm{TL}}\) may not always perform well. For instance, if the source samples have too few observations or are generated adversarially from models satisfying \(\alpha_\delta < \alpha_m\), we cannot gain any benefits from using them. In such cases, the conventional learning algorithm \(\Acal_{\mathrm{CL}}\) can still be used in the transfer learning setup by simply ignoring the collected source samples. The optimal algorithm for conventional learning, \(\Acal_{\mathrm{CL}}\), thus provides a baseline for estimation performance, as it has been shown to achieve the convergence rate of \(\widetilde{\Theta}(m_t^{-2\alpha_m}+n_t^{-1})\). The transfer learning algorithm \(\Acal_{\mathrm{TL}}\) is only advantageous if it performs better than this baseline rate; otherwise, it makes sense to stick with the conventional learning algorithm \(\Acal_{\mathrm{CL}}\). Reflecting this idea, the methodology suggested in Theorem \ref{thm:upper bound under a common design} appropriately combines \(\Acal_{\mathrm{TL}}\) with \(\Acal_{\mathrm{CL}}\) to yield an optimal estimate.

\begin{theorem}[Upper bound under a common design] \label{thm:upper bound under a common design}
Suppose
\[ \max_{1\leq j\leq m_t+1} \bigl(T_{j}^{[t]}-T_{j-1}^{[t]}\bigr) \leq \frac{C_t}{m_t} \qquad\text{and}\qquad \max_{1\leq j\leq m_s+1} \bigl(T_{j}^{[s]}-T_{j-1}^{[s]}\bigr) \leq \frac{C_s}{m_s} \]
where \(C_t, C_s > 0\) are some constants as well as \(T_0^{[t]} = T_{0}^{[s]} = 0\) and \(T_{m_t+1}^{[t]} = T_{m_t+1}^{[s]} = 1\). Consider the conventional learning estimator \(\widehat{f}_{\mathrm{CL}}^{[t]}\) which is the output of the algorithm \(\Acal_{\mathrm{CL}}(b_t, d_t, M_t)\) with the following specifications:
\begin{itemize}
    \item Any constant degree \(d_t \geq \omega(\alpha_m)\), 
    \item Any bandwidth \(b_t = \ceilbig{m_t/2B_t(d_t+1)}^{-1}\) for constant \(B_t \geq C_t\),
    \item The threshold \(M_t = \log n_t\).
\end{itemize}
Plus, let \(\widehat{f}_{\mathrm{TL}}^{[t]}\) be the output of the transfer learning algorithm \(\Acal_{\mathrm{TL}}(b_s, b_\delta, d_s, d_\delta, M_s, M_\delta)\) with the following specifications:
\begin{itemize}
    \item Any constant degrees \(d_s \geq \omega(\alpha_m)\) and \(d_\delta \geq \omega(\alpha_\delta)\).
    \item Any bandwidths \(b_s = \ceilbig{m_s/2(d_s+1)B_s}^{-1}\) and \(b_\delta = \ceilbig{m_t/2(d_\delta+1)B_\delta}^{-1}\) for any given constants \(B_s \geq C_s\) and \(B_\delta \geq C_t\).
    \item Two thresholds \(M_s = \log n_s\) and \(M_\delta = \log n_t n_s\).
\end{itemize}
The estimator \(\widehat{f}^{[t]}\) for the mean function \(f^{[t]}\) is now defined as one of them:
\[\widehat{f}^{[t]} = \begin{dcases}
    \widehat{f}_{\mathrm{CL}}^{[t]} &\text{when}~ R_{\mathrm{C}}(\widehat{f}_{\mathrm{CL}}^{[t]}) \leq R_{\mathrm{C}}(\widehat{f}_{\mathrm{TL}}^{[t]}), \\
    \widehat{f}_{\mathrm{TL}}^{[t]} &\text{when}~ R_{\mathrm{C}}(\widehat{f}_{\mathrm{CL}}^{[t]}) > R_{\mathrm{C}}(\widehat{f}_{\mathrm{TL}}^{[t]}),
\end{dcases}\]
where the following quantities are additionally defined:
\begin{align*}
    R_{\mathrm{C}}(\widehat{f}_{\mathrm{CL}}^{[t]}) &:= L_m^2 m_t^{-2\alpha_m} + \frac{\log^2 n_t}{n_t}, \\
    R_{\mathrm{C}}(\widehat{f}_{\mathrm{TL}}^{[t]}) &:= L_\delta^2 m_t^{-2\alpha_\delta} + \frac{\log^2 (n_t n_s)}{n_t} + L_m^2 m_s^{-2\alpha_m} + \frac{\log^2 (K n_s)}{K n_s}.
\end{align*}
In this situation, we obtain
\[ \sup_{\Pbb\in\Pcal} \Ebb\norm{\widehat{f}^{[t]} - f^{[t]}}^2_{\Lcal^2} \lesssim R_{\mathrm{C}}(\widehat{f}_{\mathrm{CL}}^{[t]}) \binmin R_{\mathrm{C}}(\widehat{f}_{\mathrm{TL}}^{[t]}). \]
\end{theorem}

Theorem \ref{thm:upper bound under a common design} highlights a crucial aspect of the transfer learning algorithm \(\Acal_{\mathrm{TL}}\), which achieves the convergence rate of \(\widetilde{\Theta}(m_t^{-2\alpha_\delta}+n_t^{-1}+m_s^{-2\alpha_m}+(Kn_s)^{-1})\) with appropriately chosen parameters. Corresponding to the decomposition \(f^{[t]} = f^{[s]} + \delta^{[s]}\) of the target mean function, this rate can also be decomposed into two parts: the optimal estimation error of the source part \(\widetilde{\Theta}(m_s^{-2\alpha_m}+(K n_s)^{-1})\) and the difference part \(\widetilde{\Theta}(m_t^{-2\alpha_\delta}+n_t^{-1})\). The final estimator \(\widehat{f}^{[t]}\) is now selected as an output of any algorithm that achieves a better convergence rate. We therefore obtain the best possible rate of convergence utilizing both learning algorithms, \(\Acal_{\mathrm{TL}}\) and \(\Acal_{\mathrm{CL}}\).

It should be noted that both conventional learning estimator \(\widehat{f}_{\mathrm{CL}}^{[t]}\) and the transfer learning estimator \(\widehat{f}_{\mathrm{TL}}^{[t]}\) never depend on unknown parameters. This means that both algorithms \(\Acal_{\mathrm{TL}}\) and \(\Acal_{\mathrm{CL}}\) are naturally adaptive under a common design with optimally chosen parameters. However, selecting between two outputs from two algorithms is still challenging, and the strategy for the final estimator \(\widehat{f}^{[t]}\) suggested in Theorem \ref{thm:upper bound under a common design} is not effective. In practice, the smoothness parameters \(\alpha_m\) and \(\alpha_\delta\) are typically unknown, which makes it impossible to compare the upper bounds of maximal risk, \(R_{\mathrm{C}}(\widehat{f}_{\mathrm{TL}}^{[t]})\) and \(R_{\mathrm{C}}(\widehat{f}_{\mathrm{CL}}^{[t]})\). Therefore, another adaptive and optimal procedure is necessary to decide whether to employ source samples and transfer learning. One novel solution to this issue will be presented in Section \ref{sec:adaptive estimation under a common design}.

\subsection{Matching lower bound} \label{sec:lower bound under a common design}

The following theorem demonstrates a lower bound of the maximal risk in estimating the target mean function under a common design.

\begin{theorem}[Lower bound under a common design]\label{thm:lower bound under a common design}
    Under a common design,
    \[ \inf_{\widehat{f}^{[t]}} \sup_{\Pbb\in\Pcal} \Ebb\norm{\widehat{f}^{[t]} - f^{[t]}}^2_{\Lcal^2} \gtrsim \left(L_m^2 m_t^{-2\alpha_m} + \frac{1}{n_t} \right) \binmin \left(L_\delta^2 m_t^{-2\alpha_\delta} + \frac{1}{n_t} + L_m^2 m_s^{-2\alpha_m} + \frac{1}{K n_s}\right), \]
    where the infimum is taken over all possible estimators \(\widehat{f}^{[t]} = \widehat{f}^{[t]}(\Dcal^{[t]}, \Dcal^{[s,1]}, \ldots, \Dcal^{[s,K]})\).
\end{theorem}

The lower bound presented in Theorem \ref{thm:lower bound under a common design} matches exactly with the upper bound shown in Theorem \ref{thm:upper bound under a common design}, up to logarithmic terms. This establishes the minimax risk under a common design given by:
\[\left(m_t^{-2\alpha_m} + \frac{1}{n_t} \right) \binmin \left(m_t^{-2\alpha_\delta} + \frac{1}{n_t} + m_s^{-2\alpha_m} + \frac{1}{K n_s}\right)\]
This rate is decomposed into two terms, each quantifying the difficulty of estimating the target mean function \(f^{[t]}\) but with different scenarios:
\begin{itemize}
    \item \(\widetilde{\Theta}\bigl(m_t^{-2\alpha_m}+n_t^{-1}\bigr)\) represents the conventional learning setup.
    \item \(\widetilde{\Theta}\bigl(m_t^{-2\alpha_\delta}+n_t^{-1}+m_s^{-2\alpha_m}+(K n_s)^{-1}\bigr)\) stands for the transfer learning setup.
\end{itemize}
It is reasonable to measure the entire difficulty of the problem by the minimum of these two terms. Furthermore, the second term can be decomposed into \(\widetilde{\Theta}\bigl(m_s^{-2\alpha_m}+(Kn_s)^{-1}\bigr)\) and \(\widetilde{\Theta}\bigl(m_t^{-2\alpha_\delta}+n_t^{-1}\bigr)\), which correspond to the minimax risk of estimation on the source and target samples, respectively. This decomposition highlights the crux of transfer learning: by paying the cost of estimation on the source sample, the target sample can be treated as if the target mean function is given by \(\alpha_\delta\)-smooth, instead of \(\alpha_m\)-smooth. As long as the cost is not exceedingly expensive and substituting smoothness results in increased smoothness, this leads to a more accurate estimation and a faster rate of convergence.

By comparing the minimax risks between the conventional and transfer learning settings, we can determine whether and under what conditions the source samples and transfer learning can enhance performance. The comparison is relatively simple because the minimax rate for the transfer learning setup already includes that for the conventional learning setup. What makes this process particularly interesting is the presence of a phase transition under which the efficacy of transfer learning experiences a significant shift.

To begin with, there is no rationale for incurring the cost of substituting smoothness if it does not result in increased smoothness. In other words, if our model is given with \(\alpha_\delta \leq \alpha_m\), transfer learning has an adverse effect because estimating \(\alpha_\delta\)-smooth functions is, at the very least, as challenging as estimating \(\alpha_m\)-smooth functions. The effectiveness of transfer learning is justifiable only for models such that \(\alpha_\delta\) is strictly greater than \(\alpha_m\). In such cases, we can exploit the benefits of substituting smoothness and anticipate surpassing the minimax risk of the conventional learning algorithm. This gives rise to the phase transition phenomenon and further discussions about the effectiveness of transfer learning become meaningful under the additional assumption that \(\alpha_\delta > \alpha_m\) holds.

From this point onwards, suppose that our model satisfies \(\alpha_\delta > \alpha_m\). In the high-frequency regime, defined as \(m_t \gtrsim n_t^{1/2\alpha_m}\), the minimax risk for the conventional learning algorithm is \(\widetilde{\Theta}(n_t^{-1})\) because this rate is always as rapid as \(\widetilde{\Theta}(m_t^{-2\alpha_\delta}+n_t^{-1}+m_s^{-2\alpha_m}+(Kn_s)^{-1})\), the additional rate introduced by the transfer learning setup. In other words, the conventional learning algorithm already achieves the parametric rate \(\widetilde{\Theta}(n_t^{-1})\) concerning the number of target subjects, regardless of the existence of source samples or the level of smoothness. Unless difference functions \(\delta^{[s,k]}\) \((k=1,\ldots, K)\) are exactly known, the minimax risk for estimating the target mean function must be at least the parametric rate \(\widetilde{\Theta}(n_t^{-1})\). As a result, transfer learning cannot be effective as there is no more room for improvement in the high-frequency regime. This kind of phenomenon is expected to be universal across problems beyond the scope of the functional mean estimation.

In the low-frequency regime with \(m_t \ll n_t^{1/2\alpha_m}\), the minimax risk in the conventional learning setting is \(\widetilde{\Theta}(m_t^{-2\alpha_m})\). In this case, transfer learning improves the estimation performance of the target mean function provided the following condition holds:
\[\widetilde{\Theta}\bigl(m_s^{-2\alpha_m} + (Kn_s)^{-1}\bigr) \ll \widetilde{\Theta}\bigl(m_t^{-2\alpha_m}\bigr).\]
This condition implies that estimating the average of source mean functions denoted as \(f^{[s]}\), is a less challenging task compared to estimating the target mean function \(f^{[t]}\) assuming the conventional learning settings. Moreover, this condition remains unaffected by the smoothness parameter \(\alpha_\delta\), indicating that the difference functions \(\delta^{[s,k]}\) \((k=1,\ldots, K)\) only need to be slightly smoother than the target mean function \(f^{[t]}\).

In breaking down the condition for the low-frequency regime, we can identify two essential requirements. First, the source subjects in total should be abundant, satisfying \(K n_s \gg m_t^{2\alpha_m}\). Second, the frequency of design points in the source samples, \(m_s\), should significantly surpass that of the target sample, \(m_t\), indicating \(m_s \gg m_t\). It is worth noting that the size of the source group denoted as \(K\), plays a restricted role under a common design. In cases where the source samples lack an adequate number of design points \((m_s \lesssim m_t)\), simply increasing the size of the source group will not yield substantial benefits in transfer learning. As long as the source samples offer a sufficient pool of design points \((m_s \gg m_t)\), increasing \(K\) can be an effective strategy to obtain a faster rate of convergence with transfer learning.

\subsection{Adaptive estimation} \label{sec:adaptive estimation under a common design}

Although the optimal estimators \(\widehat{f}_{\mathrm{CL}}^{[t]}\) and \(\widehat{f}_{\mathrm{TL}}^{[t]}\) introduced in Theorem \ref{thm:upper bound under a common design} are naturally adaptive, selecting between them is not. We propose a new adaptive procedure to address this issue. Algorithm \ref{alg:adaptive learning algorithm under a common design} is an outline of the adaptive learning algorithm under a common design, called \(\Acal_{\mathrm{ALC}}\). It is further assumed for brevity that the target sample \(\Dcal^{[t]}\) contains \(2n_t\) subjects, which does not affect the rate of convergence.

\begin{algorithm}[H]
    \caption{Adaptive transfer learning for mean function under a common design \(\Acal_{\mathrm{ALC}}\)} \label{alg:adaptive learning algorithm under a common design}
    \begin{algorithmic}[1] \normalsize
    \State Randomly partition the target sample \(\Dcal^{[t]}\) into two sub-samples, denoted as \(\Dcal_{\mathrm{train}}^{[t]}\) and \(\Dcal_{\mathrm{test}}^{[t]}\), based on subjects. Specifically, perform a random split of the index set \(\{1,\ldots,2n_t\}\) into two partitions, \(\Ical_{\mathrm{train}}^{[t]}\) and \(\Ical_{\mathrm{test}}^{[t]}\) such that \(\abs{\Ical_{\mathrm{train}}^{[t]}} = \abs{\Ical_{\mathrm{test}}^{[t]}} = n_t\). We define:
    \[ \begin{aligned}
        \Dcal_{\mathrm{train}}^{[t]} &:= \bigl\{(T_{j}^{[t]}, Y_{ij}^{[t]}): i \in \Ical_{\mathrm{train}}^{[t]},~ j = 1,\ldots,m_t\bigr\}, \\
        \Dcal_{\mathrm{test}}^{[t]} &:= \bigl\{(T_{j}^{[t]}, Y_{ij}^{[t]}): i \in \Ical_{\mathrm{test}}^{[t]},~ j = 1,\ldots,m_t\bigr\}.
    \end{aligned} \]
    \State Execute both \(\Acal_{\mathrm{CL}}(b_t, d_t, M_t)\) and \(\Acal_{\mathrm{TL}}(b_s, b_\delta, d_s, d_\delta, M_s, M_\delta)\) following the same specifications as outlined in Theorem \ref{thm:upper bound under a common design}, with the only distinction being that the target sample is provided as \(\Dcal_{\mathrm{train}}^{[t]}\), not \(\Dcal^{[t]}\). The outputs are denoted as \(\widehat{f}_{\mathrm{CL}}^{[t]}\) and \(\widehat{f}_{\mathrm{TL}}^{[t]}\), respectively.
    \State Output the following estimator \(\widehat{g}_*^{[t]}\). If a tie occurs, use any randomization to break it.
    \[ \widehat{g}_*^{[t]} := \argmin_{\widehat{g}^{[t]} \in \{\widehat{f}_{\mathrm{CL}}^{[t]},\widehat{f}_{\mathrm{TL}}^{[t]}\}} \sum_{i \in \Ical_{\mathrm{test}}^{[t]}} \sum_{j=1}^{m_t} \bigl(Y^{[t]}_{ij} - \widehat g^{[t]}(T_{j}^{[t]})\bigr)^2 (\Delta T_j^{[t]}), \]
    where \(T_0^{[t]} := 0\), \(T_{m_t+1}^{[t]} := 1\) and \(\Delta T_j^{[t]} := T_{j}^{[t]} - T_{j-1}^{[t]}\) for each \(j=1,\ldots,m_t+1\).
\end{algorithmic}
\end{algorithm}

The success of Algorithm \ref{alg:adaptive learning algorithm under a common design} mainly hinges on the train-test split. Half of the target sample, denoted as \(\Dcal_{\mathrm{train}}^{[t]}\), serves as the training data for two algorithms, while the remaining portion, referred to as \(\Dcal_{\mathrm{test}}^{[t]}\), is employed to evaluate their empirical performance. Although the design points of the target sample are fixed, we have assumed \(\max_{j=1,\ldots,m_t+1} (\Delta T_j^{[t]}) = o(m_t^{-1})\), ensuring that the Riemann sum serves as a robust proxy for the functional \(\Lcal^2\)-loss. Not only that, the computational cost for the adaptation step remains negligible, as it merely involves the selection between two data-driven estimators, \(\widehat{f}_{\mathrm{CL}}^{[t]}\) and \(\widehat{f}_{\mathrm{TL}}^{[t]}\).

Notice that the base algorithms \(\Acal_{\mathrm{TL}}\) and \(\Acal_{\mathrm{CL}}\) generates a randomized estimator. Consequently, the adaptive learning algorithm \(\Acal_{\mathrm{ALC}}\) must be randomized as well. To improve the performance in finite samples, we execute this adaptive procedure \(r_{\mathrm{max}}\) times and calculate the average of the resulting estimates. This approach shares a similar idea with bagging estimators suggested by \citet{breimanBaggingPredictors1996}. The optimal minimax rate of convergence, up to logarithmic terms, is achieved by this algorithm, as demonstrated by the subsequent theorem.

\begin{theorem}[Adaptive estimation under a common design] \label{thm:adaptive estimation under a common design}
    Under the same assumptions as Theorem \ref{thm:upper bound under a common design}, we consider a maximum number of repetitions, \(r_{\mathrm{max}}\in\Zbb^+\) and let \(\widehat{g}_r^{[t]}\) \((r=1,\ldots,r_{\mathrm{max}})\) denote the output of the \(r\)-th execution of the algorithm \(\Acal_{\mathrm{ALC}}\). By averaging these estimates, we obtain the final estimator:
    \[ \widehat{f}^{[t]} = \frac{1}{r_{\mathrm{max}}} \sum_{r=1}^{r_{\mathrm{max}}} \widehat{g}_r^{[t]}. \numberthis\label{eq:adaptive common design} \]
    This adaptive estimator \(\widehat{f}^{[t]}\) attains the same upper bound of Theorem \ref{thm:upper bound under a common design}:
    \[ \sup_{\Pbb\in\Pcal} \Ebb\norm{\widehat{f}^{[t]} - f^{[t]}}^2_{\Lcal^2} \lesssim R_{\mathrm{C}}(\widehat{f}_{\mathrm{CL}}^{[t]}) \binmin R_{\mathrm{C}}(\widehat{f}_{\mathrm{TL}}^{[t]}). \]
\end{theorem}
In short, the data-driven estimator \eqref{eq:adaptive common design} adaptively achieves the optimal rate of convergence over a comprehensive collection of function classes.

%% file: content/independent-design.tex

\section{Transfer learning for functional mean estimation under an independent design} \label{sec:optimality under an independent design}

This section investigates the setting in which the design points of both the target and source samples are randomly drawn from the domain \([0,1]\). Consider two probability distributions on \([0,1]\), referred to as \(\eta_t\) and \(\eta_s\), for target and sample design points, respectively. We make the following assumptions:
\[ \begin{aligned}
    &\bigl(T_{ij}^{[t]}: j=1,\ldots,m_t,~ i=1,\ldots,n_t \bigr) \iid \eta_t, \\
    &\bigl(T_{ij}^{[s,k]}: j=1,\ldots,m_s,~ i=1,\ldots,n_s,~ k=1,\ldots,K \bigr) \iid \eta_s. \\
\end{aligned} \]

Apart from the assumption regarding design points, we maintain the same set of assumptions outlined in Section \ref{sec:optimality under a common design}. In essence, these assumptions encompass the requirements for the mean and difference functions, as indicated in Equation \eqref{eq:modeling assumption for mean and difference functions}, and the fulfillment of the uniform sub-Gaussian condition (Assumption \ref{ass:uniformly subgaussian assumption}). Lastly, it is assumed that the collection of target and source samples, \(\{\Dcal^{[t]}, \Dcal^{[s,1]}, \ldots, \Dcal^{[s,K]}\}\), are independent. With these assumptions at hand, we can proceed to examine the optimal estimation of the target mean function \(f^{[t]}\) under an independent design.

\subsection{Methodology and upper bound}

This subsection introduces an optimal algorithm for estimating the target mean function. In parallel to the discussion for a common design, let us reconsider the conventional setting for functional mean estimation but this time under an independent design. We again aim to estimate the target mean function \(f^{[t]}\) when we have no source samples available \((n_s = 0)\). The subsequent theorem investigates the minimax rate of convergence in this situation.

\begin{theorem}[The minimax risk under conventional setup and independent design] \label{thm:minimax risk under conventional setup and independent design}
    Suppose the distribution \(\eta_t\) for the target design points is dominated by the Lebesgue measure and its density is bounded from below by a constant \(C_t > 0\) and from above by \(\Gamma_t > 0\). Then
    \[ \inf_{\widehat{f}^{[t]}} \sup_{\Pbb\in\Pcal} \Ebb\norm{\widehat{f}^{[t]} - f^{[t]}}^2_{\Lcal^2} = \widetilde{\Theta} \left(L_m^{2/(2\alpha_m+1)} (m_t n_t)^{-2\alpha_m/(2\alpha_m+1)} + \frac{1}{n_t}\right), \]
    where the infimum is taken over all estimators \(\widehat{f}^{[t]} = \widehat{f}^{[t]}(\Dcal^{[t]})\) based on the target sample.
\end{theorem}

The derived minimax rate of convergence again generalizes the optimal rate shown by \citet{caiOptimalEstimationMean2011}, \(\Theta((m_t n_t)^{-2r/(2r+1)}+n_t^{-1})\), where they assumed that both the curve \(X^{[t]}\) and the mean function \(f^{[t]}\) are differentiable \(r\) times \((r\in\Zbb^+)\). Our framework provides greater flexibility and improvement over previous approaches for functional mean estimation, similar to what we have observed under a common design.

Even though the sampling designs are different, the current model has a similar structure to the corresponding one under a common design. We thus employ the same conventional learning algorithm \(\Acal_{\mathrm{CL}}\) for estimating the target mean function. The output of the algorithm, denoted as \(\widehat{f}_{\mathrm{CL}}^{[t]}\), achieves the minimax rate of convergence in Theorem \ref{thm:minimax risk under conventional setup and independent design}, but the optimal choice of input parameters is distinct from that of a common design. We will provide more details about this choice in Theorem \ref{thm:upper bound under an independent design}.
\[ \sup_{\Pbb\in\Pcal} \Ebb\norm{\widehat{f}_{\mathrm{CL}}^{[t]} - f^{[t]}}_{\Lcal^2}^2 \leq O\left(L_m^{2/(2\alpha_m+1)} (m_t n_t)^{-2\alpha_m/(2\alpha_m+1)} + \frac{\log^2 n_t}{n_t}\right). \]

We will next shift our attention back to the transfer learning setup of our primary focus. Our approach for estimating the target mean function involves implementing and combining the two learning algorithms \(\Acal_{\mathrm{CL}}\) and \(\Acal_{\mathrm{TL}}\) optimally. The underlying idea is that the transfer learning algorithm \(\Acal_{\mathrm{TL}}\) may not always perform well, especially when the source samples are limited in number or generated in an adversarial way. In such cases, we can exploit the conventional learning algorithm \(\Acal_{\mathrm{CL}}\) as a baseline. This approach is again not surprising as an analogous type of estimator has been demonstrated to achieve the matching lower bound under a common design. The optimal way to execute and combine both learning algorithms under an independent design is outlined in Theorem \ref{thm:upper bound under an independent design}. Although the overall structure of the estimator is comparable to that of a common design, the optimal selection of parameters is substantially different.

\begin{theorem}[Upper bound under an independent design] \label{thm:upper bound under an independent design}
    Suppose the distributions \(\eta_t\) and \(\eta_s\) for target and source design points are dominated by the Lebesgue measure, and their densities are bounded from below by a constant \(C_t, C_s > 0\) and from above by \(\Gamma_t, \Gamma_s > 0\), respectively. Consider the conventional learning estimator \(\widehat{f}_{\mathrm{CL}}^{[t]}\) which is the output of algorithm \(\Acal_{\mathrm{CL}}(b_t, d_t, M_t)\) with the following specifications:
    \begin{itemize}
        \item For any constant \(B_t \geq \Gamma_t\), bandwidth:
        \[ b_t = \ceilauto{(L_m^2 m_t n_t)^{1/(2\alpha_m + 1)} (\log n_t)^{-2/(2\alpha_m+1)}}^{-1} \binmin \ceilbig{2B_t m_t}^{-1}, \]
        \item Any constant degree \(d_t \geq \omega(\alpha_m)\),
        \item The threshold \(M_t = \log n_t\).
    \end{itemize}
    Besides, let \(\widehat{f}_{\mathrm{TL}}^{[t]}\) be the output of the transfer learning algorithm \(\Acal_{\mathrm{TL}}(b_s, b_\delta, d_s, d_\delta, M_s, M_\delta)\) with the following specifications:
    \begin{itemize}
        \item For any constants \(B_s \geq \Gamma_s\) and \(B_\delta \geq \Gamma_t\), bandwidths:
        \begin{align*}
            b_s &= \ceilauto{(L_m^2 K m_s n_s)^{1/(2\alpha_m + 1)} \bigl(\log (K n_s)\bigr)^{-2/(2\alpha_m+1)}}^{-1} \binmin \ceilbig{2B_s K m_s}^{-1}, \\
            b_\delta &= \ceilauto{(L_\delta^2 m_t n_t)^{1/(2\alpha_\delta + 1)} \bigl(\log (n_t n_s)\bigr)^{-2/(2\alpha_\delta+1)}}^{-1} \binmin \ceilbig{2B_\delta m_t}^{-1}.
        \end{align*}
        \item Any constant degrees \(d_s \geq \omega(\alpha_m)\) and \(d_\delta \geq \omega(\alpha_\delta)\).
        \item Two thresholds \(M_s = \log n_s\) and \(M_\delta = \log n_t n_s\).
    \end{itemize}
    The estimator \(\widehat{f}^{[t]}\) of the mean function \(f^{[t]}\) is defined as one of them:
    \[\widehat{f}^{[t]} = \begin{dcases}
        \widehat{f}_{\mathrm{CL}}^{[t]} &\text{when}~ R_{\mathrm{I}}(\widehat{f}_{\mathrm{CL}}^{[t]}) \leq R_{\mathrm{I}}(\widehat{f}_{\mathrm{TL}}^{[t]}), \\
        \widehat{f}_{\mathrm{TL}}^{[t]} &\text{when}~ R_{\mathrm{I}}(\widehat{f}_{\mathrm{CL}}^{[t]}) > R_{\mathrm{I}}(\widehat{f}_{\mathrm{TL}}^{[t]}),
    \end{dcases}\]
    where the following quantities are further defined:
    \begin{align*}
        R_{\mathrm{I}}(\widehat{f}_{\mathrm{CL}}^{[t]}) &:= L_m^{2/(2\alpha_m+1)} \left(\frac{\log^2 n_t}{m_t n_t}\right)^{\!2\alpha_m/(2\alpha_m+1)} + \frac{\log^2 n_t}{n_t}, \\
        R_{\mathrm{I}}(\widehat{f}_{\mathrm{TL}}^{[t]}) &:= L_m^{2/(2\alpha_m+1)} \left(\frac{\log^2 (K n_s)}{K m_s n_s}\right)^{\!2\alpha_m/(2\alpha_m+1)} + \frac{\log^2 (K n_s)}{K n_s} \\
        &~~~~~~+ L_\delta^{2/(2\alpha_\delta+1)} \left(\frac{\log^2 (n_t n_s)}{m_t n_t}\right)^{\!2\alpha_\delta/(2\alpha_\delta+1)} + \frac{\log^2 (n_t n_s)}{n_t}.
    \end{align*}
    In this situation, we have
    \[ \sup_{\Pbb\in\Pcal} \Ebb\norm{\widehat{f}^{[t]} - f^{[t]}}^2_{\Lcal^2} \lesssim R_{\mathrm{I}}(\widehat{f}_{\mathrm{CL}}^{[t]}) \binmin R_{\mathrm{I}}(\widehat{f}_{\mathrm{TL}}^{[t]}). \]
\end{theorem}

It is worth noting in Theorem \ref{thm:upper bound under an independent design} that the transfer learning algorithm \(\Acal_{\mathrm{TL}}\) achieves the rate \(\widetilde{\Theta}\bigl((m_t n_t)^{-2\alpha_\delta/(2\alpha_\delta+1)} + n_t^{-1}+(K m_s n_s)^{-2\alpha_m/(2\alpha_m+1)}+(K n_s)^{-1}\bigr)\), provided that the input parameters are selected appropriately. However, it is not always superior to the conventional learning algorithm, which has an optimal rate of \(\widetilde{\Theta}\bigl((m_t n_t)^{-2\alpha_m/(2\alpha_m+1)}+n_t^{-1}\bigr)\). Therefore, the final estimator \(\widehat{f}^{[t]}\) in Theorem \ref{thm:upper bound under an independent design} is selected as any algorithm that achieves a better rate of convergence, which optimizes the utilization of both learning algorithms \(\Acal_{\mathrm{TL}}\) and \(\Acal_{\mathrm{CL}}\).

Although the final estimator \(\widehat{f}^{[t]}\) suggested in Theorem \ref{thm:upper bound under an independent design} will be sufficient to argue the minimax rate of convergence, it confronts two challenges in real-world situations. Since the smoothness parameters \(\alpha_m\) and \(\alpha_\delta\) are typically unknown, it is impossible to obtain the optimal bandwidth for both algorithms in practice. Additionally, comparing the rates of the two algorithms also requires knowledge of those unknown parameters. Adaptation to smoothness parameters was relatively simple under a common design since the optimal parameters for both learning algorithms are at least independent of them. As a consequence, it must be more challenging to decide whether or not to employ the source samples and transfer learning under an independent design.

\subsection{Matching lower bound}

The upcoming theorem shows the lower bound for estimating the target mean function under an independent design.

\begin{theorem}[Lower bound under an independent design]\label{thm:lower bound under an independent design}
    Under an independent design,
    \begin{align*}
        \inf_{\widehat{f}^{[t]}} \sup_{\Pbb\in\Pcal} \Ebb\norm{\widehat{f}^{[t]} - f^{[t]}}^2_{\Lcal^2} &\gtrsim \left[L_m^{2/(2\alpha_m+1)} (m_t n_t)^{-2\alpha_m/(2\alpha_m+1)} + \frac{1}{n_t}\right] \\
        &~~~~~~\binmin \biggl[ L_\delta^{2/(2\alpha_\delta+1)}(m_t n_t)^{-2\alpha_\delta/(2\alpha_\delta+1)} + \frac{1}{n_t} \biggr. \\
        &~~~~~~~~~~~~~~\biggl. + L_m^{2/(2\alpha_m+1)}(K m_s n_s)^{-2\alpha_m/(2\alpha_m+1)} + \frac{1}{K n_s} \biggr],
    \end{align*}
    where the infimum is taken over all possible estimators \(\widehat{f}^{[t]} = \widehat{f}^{[t]}(\Dcal^{[t]}, \Dcal^{[s,1]}, \ldots, \Dcal^{[s,\ell]})\).
\end{theorem}

Theorem \ref{thm:lower bound under an independent design} provides the minimax lower bound which exactly matches the upper bound given in Theorem \ref{thm:upper bound under an independent design}, up to logarithmic factors. This leads to the optimal rate for estimating the target mean function \(f^{[t]}\) under an independent design, which is given by:
\[ \left[(m_t n_t)^{-2\alpha_m/(2\alpha_m+1)} + \frac{1}{n_t}\right] \binmin \left[(m_t n_t)^{-2\alpha_\delta/(2\alpha_\delta+1)} + \frac{1}{n_t} + (K m_s n_s)^{-2\alpha_m/(2\alpha_m+1)} + \frac{1}{K n_s} \right] \]
This rate comprises two primary components:
\begin{itemize}
    \item \(\widetilde{\Theta}\bigl((m_t n_t)^{-2\alpha_m/(2\alpha_m+1)}+n_t^{-1}\bigr)\), arising from the conventional learning setup.
    \item \(\widetilde{\Theta}\bigl((m_t n_t)^{-2\alpha_\delta/(2\alpha_\delta+1)}+n_t^{-1}+(K m_s n_s)^{-2\alpha_m/(2\alpha_m+1)} + (Kn_s)^{-1}\bigr)\), originating from the transfer learning setup.
\end{itemize}
It makes sense to evaluate the problem's overall difficulty by taking the minimum of these two rates. The latter rate is further decomposed into two parts: \(\widetilde{\Theta}\bigl((m_t n_t)^{-2\alpha_\delta/(2\alpha_\delta+1)}+n_t^{-1}\bigr)\) and \(\widetilde{\Theta}\bigl((K m_s n_s)^{-2\alpha_m/(2\alpha_m+1)}+(Kn_s)^{-1}\bigr)\). These rates correspond to the minimax risks of estimation on the target and source samples, respectively. In parallel to the discussion under a common design, transfer learning treats the target mean function as \(\alpha_\delta\)-smooth, instead of \(\alpha_m\)-smooth, while paying the cost of estimating the source mean functions. Whenever this cost is reasonable and the substitution of smoothness results in benefits, it leads to more accurate estimation and a faster convergence rate. This high-level structure of the minimax risk is similar to that of a common design and reveals that the model's formulation does not depend on the sampling designs.

By analyzing the minimax risks within the contexts of both the conventional and transfer learning settings, we can determine the circumstances under which utilizing source samples or transfer learning can enhance the performance of estimating the target mean function. It is straightforward to compare the two settings as the minimax rate for the conventional learning setup comprises half of that for the transfer learning setup. Similar to a common design case, we can identify a phase transition that has a substantial impact on the effectiveness of transfer learning. Remarkably, the boundary of phase transition is the same as that in a common design setting. This also implies that the same fundamental principles govern the functional mean estimation regardless of the sampling designs.

First and foremost, it suffices to focus our analysis of effective transfer learning exclusively on models that adhere to the condition \(\alpha_\delta > \alpha_m\). The rationale behind this criterion aligns with that of the common design setting. Within such models, we may fully exploit the benefits of substituting smoothness; the increase in smoothness potentially leads to a faster rate of convergence than the one associated with the conventional learning algorithm. For any model with \(\alpha_\delta \leq \alpha_m\), transfer learning indeed yields adverse effects because the estimation of \(\alpha_\delta\)-smooth functions is as challenging as estimating \(\alpha_m\)-smooth functions. Under the absence of clear advantage in smoothness, there is no justification for engaging in estimation on the source samples and incurring the cost of \(\widetilde{\Theta}\bigl((K m_s n_s)^{-2\alpha_m/(2\alpha_m+1)}+(Kn_s)^{-1}\bigr)\) in the convergence rate. This introduces the intriguing concept of a phase transition where the effectiveness of transfer learning undergoes a sharp change. Further discussions concerning the effectiveness become meaningful for limited models satisfying \(\alpha_\delta > \alpha_m\).

We exclusively focus on models that meet the condition \(\alpha_\delta > \alpha_m\). In the high-frequency regime where \(m_t \gtrsim n_t^{1/2\alpha_m}\), transfer learning cannot be effective, similar to the case of a common design. The convergence rate for the conventional learning setup, \(\widetilde{\Theta}(n_t^{-1})\), is always less than or equal to \(\widetilde{\Theta}\bigl((m_t n_t)^{-2\alpha_\delta/(2\alpha_\delta+1)}+n_t^{-1}+(K m_s n_s)^{-2\alpha_m/(2\alpha_m+1)}+(K n_s)^{-1}\bigr)\), the additional rate introduced by the transfer learning setup. Unless difference functions \(\delta^{[s,k]}\) \((k=1,\ldots,\ell)\) are exactly known, the parametric rate \(\widetilde{\Theta}(n_t^{-1})\) always gives the best possible rate and there is no further way for enhancement.

However, in the low-frequency regime where \(m_t \ll n_t^{1/2\alpha_m}\), transfer learning has the potential to enhance the estimation performance like the case of a common design. To be more specific, the following condition is necessary and sufficient for the source samples and transfer learning to be beneficial:
\[ (m_s n_s)^{-2\alpha_m/(2\alpha_m+1)}+(Kn_s)^{-1} \ll (m_t n_t)^{-2\alpha_m/(2\alpha_m+1)}. \]
This condition is quite different from the corresponding one under a common design, which is \(m_s^{-2\alpha_m}+n_s^{-1}\ll m_t^{-2\alpha_m}\), and neither condition implies the other. Comparing these two conditions for effective transfer learning, both require a sufficient number of source subjects, but their thresholds are different. For an independent design, the requirement is related to the target observation size \((K n_s \gg (m_t n_t)^{2\alpha_m/(2\alpha_m+1)})\), while for a common design, it is to the target sampling frequency \((K n_s \gg m_t^{2\alpha_m})\). Additionally, an independent design requires the source samples to have more total observations \((K m_s n_s \gg m_t n_t)\), whereas a common design requires the source samples to have a higher sampling frequency \((m_s \gg m_t)\). 

The impact of the source group's size \(K\) varies between a common design and an independent design. As previously discussed, in a common design, increasing \(K\) may or may not result in a more effective transfer learning in the low-frequency regime. However, under an independent design, augmenting \(K\) consistently leads to enhanced transfer learning in the same regime. Given that the target sample remains unchanged, it should be an effective approach for both conditions, \(K n_s \gg (m_t n_t)^{2\alpha_m/(2\alpha_m+1)}\) and \(K m_s n_s \gg m_t n_t\) to be satisfied. Although the model structures are similar across sampling designs, the minimax risks, the particular conditions for effective transfer learning, and the impact of the source group's size differ significantly. This suggests that the estimation behavior is entirely distinct between the two designs.

It is also interesting to directly compare the minimax risks between common and independent designs. Assuming that all other specifications are the same except for the sampling design, the independent design's rate is always equal to or lower than the common design's, up to some constants. Hence, if given the option to choose between the two designs, and all other factors remain the same, then an independent design should be preferred in terms of the convergence rate.
\begin{align*}
    &\left[(m_t n_t)^{-2\alpha_m/(2\alpha_m+1)} + \frac{1}{n_t}\right] \binmin \left[(m_t n_t)^{-2\alpha_\delta/(2\alpha_\delta+1)} + \frac{1}{n_t} + (K m_s n_s)^{-2\alpha_m/(2\alpha_m+1)} + \frac{1}{K n_s} \right] \\
    &~~\lesssim \left[ m_t^{-2\alpha_m} + \frac{1}{n_t} \right] \binmin \left[ m_t^{-2\alpha_\delta} + \frac{1}{n_t} + m_s^{-2\alpha_m} + \frac{1}{Kn_s} \right].
\end{align*}

\subsection{Adaptive estimation}

We introduce an adaptive algorithm for an independent design, named \(\Acal_{\mathrm{ALI}}\), which addresses two main challenges: selecting an optimal bandwidth and choosing a learning algorithm between \(\Acal_{\mathrm{CL}}\) and \(\Acal_{\mathrm{TL}}\). The algorithm generates multiple candidate bandwidths for each algorithm and compares their performances using empirical risk and cross-validation. This approach addresses both challenges simultaneously. We will provide a detailed description of \(\Acal_{\mathrm{TL}}\) in Algorithm \ref{alg:adaptive learning algorithm under an independent design}. It is further assumed for brevity that the target sample \(\Dcal^{[t]}\) contains \(2n_t\) subjects, which does not affect the rate of convergence.

It should be noted that the lists of candidate bandwidths proposed in Algorithm \ref{alg:adaptive learning algorithm under an independent design} always include the optimal bandwidths, presented in Theorem \ref{thm:upper bound under an independent design}, for the conventional and transfer learning algorithms. Therefore, any chosen bandwidth through cross-validation is expected to perform better than the optimal one. The same reasoning applies to the selection of an optimal learning algorithm. Any estimator selected via cross-validation will outperform both learning algorithms \(\Acal_{\mathrm{TL}}\) and \(\Acal_{\mathrm{CL}}\) with their respective optimal parameters. Finally, the computational complexity of the adaptation step is quite minimal, as the number of candidate bandwidths is expected to be at most \(\log_2(m_t n_t) + \log_2(Km_s n_s) \log_2(m_t n_t)\).

Similar to the adaptive learning algorithm \(\Acal_{\mathrm{ALC}}\) for a common design, the adaptive algorithm \(\Acal_{\mathrm{ALI}}\) for an independent design also involves multiple steps of random down-sampling, resulting in a randomized estimator. To improve the performance of estimation in finite samples, we will execute this adaptive process \(r_{\mathrm{max}}\) times and average the resulting estimates. This approach is analogous to bagging estimators proposed by \citet{breimanBaggingPredictors1996}. As demonstrated by the subsequent theorem, this algorithm achieves the optimal minimax rate of convergence, up to logarithmic factors.

\begin{algorithm}[H]
    \caption{Adaptive transfer learning for mean function under an independent design \(\Acal_{\mathrm{ALI}}\)} \label{alg:adaptive learning algorithm under an independent design}
    \begin{algorithmic}[1] \normalsize
    \State Initialize \(\widehat{\Gcal}^{[t]} = \emptyset\), the collection of candidate estimators.
    \State Randomly partition the target sample \(\Dcal^{[t]}\) into two sub-samples, denoted as \(\Dcal_{\mathrm{train}}^{[t]}\) and \(\Dcal_{\mathrm{test}}^{[t]}\), based on subjects. Specifically, perform a random split of the index set \(\{1,\ldots,2n_t\}\) into two partitions, \(\Ical_{\mathrm{train}}^{[t]}\) and \(\Ical_{\mathrm{test}}^{[t]}\) such that \(\abs{\Ical_{\mathrm{train}}^{[t]}} = \abs{\Ical_{\mathrm{test}}^{[t]}} = n_t\). We define:
    \[ \begin{aligned}
        \Dcal_{\mathrm{train}}^{[t]} &:= \bigl\{(T_{ij}^{[t]}, Y_{ij}^{[t]}): i \in \Ical_{\mathrm{train}}^{[t]},~ j = 1,\ldots,m_t\bigr\}, \\
        \Dcal_{\mathrm{test}}^{[t]} &:= \bigl\{(T_{ij}^{[t]}, Y_{ij}^{[t]}): i \in \Ical_{\mathrm{test}}^{[t]},~ j = 1,\ldots,m_t\bigr\}.
    \end{aligned} \]
    \State Pick any constants \(B_t \geq C_t\) and \(d_t \geq \omega(\alpha_m)\). 
    \State Take a threshold \(M_t = \log n_t\).
    \For{\(b_t \in \{2^{r} \leq m_t n_t : r \in \Zbb^+\}\)}
        \State Execute \(\Acal_{\mathrm{CL}}(b_t, d_t, M_t)\) as if target sample is given as \(\Dcal_{\mathrm{train}}^{[t]}\), not \(\Dcal^{[t]}\).
        \State Add the algorithm's output to the collection \(\widehat{\Gcal}^{[t]}\).
    \EndFor
    \State Pick any constants \(B_s \geq C_s\), \(B_\delta \geq C_t\), \(d_s \geq \omega(\alpha_m)\) and \(d_\delta \geq \omega(\alpha_\delta)\).
    \State Take thresholds \(M_s = \log n_s\) and \(M_s = \log n_t n_s\).
    \For{\((b_s, b_\delta) \in \{2^{r} \leq K m_s n_s : r \in \Zbb^+\} \times \{2^{r} \leq m_t n_t : r \in \Zbb^+\}\)}
        \State Execute \(\Acal_{\mathrm{TL}}(b_s, b_\delta, d_s, d_\delta, M_s, M_\delta)\) as if target sample is given as \(\Dcal_{\mathrm{train}}^{[t]}\), not \(\Dcal^{[t]}\).
        \State Add the algorithm's output to the collection \(\widehat{\Gcal}^{[t]}\).
    \EndFor
    \State Output the following estimator \(\widehat{g}_*^{[t]}\). If a tie occurs, use any randomization to break it.
    \[ \widehat{g}_*^{[t]} := \argmin_{\widehat{g}^{[t]} \in \widehat{\Gcal}^{[t]}} \sum_{(T,Y) \in \Dcal_{\mathrm{test}}^{[t]}} \bigl(Y - \widehat g^{[t]}(T)\bigr)^2. \]
\end{algorithmic}
\end{algorithm}

\begin{theorem}[Adaptive estimation under an independent design] \label{thm:adaptive estimation under an independent design}
    Suppose the assumptions in Theorem \ref{thm:upper bound under an independent design} hold. For any given \(r_{\mathrm{max}} \in \Zbb^+\), let \(\widehat{g}_r^{[t]}\) \((r=1,\ldots,r_{\mathrm{max}})\) be the output of the \(r\)-th execution of the algorithm \(\Acal_{\mathrm{ALI}}\). Take an average of them to obtain our final estimator:
    \[\widehat{f}^{[t]} = \frac{1}{r_{\mathrm{max}}} \sum_{r=1}^{r_{\mathrm{max}}} \widehat{g}_r^{[t]}. \numberthis\label{eq:adaptive independent design}\]
    This adaptive estimator \(\widehat{f}^{[t]}\) attains the same upper bound of Theorem \ref{thm:upper bound under an independent design}:
    \[ \sup_{\Pbb\in\Pcal} \Ebb\norm{\widehat{f}^{[t]} - f^{[t]}}^2_{\Lcal^2} \lesssim R_{\mathrm{I}}(\widehat{f}_{\mathrm{CL}}^{[t]}) \binmin R_{\mathrm{I}}(\widehat{f}_{\mathrm{TL}}^{[t]}). \]
\end{theorem}
Therefore, the data-driven estimator \eqref{eq:adaptive independent design} attains the optimal rate of convergence adaptively over a wide collection of function classes.

%% file: content/numerical-experiments.tex

\section{Numerical experiments} \label{sec:numerical experiments}

The proposed adaptive transfer learning algorithms, \(\Acal_{\mathrm{ALC}}\) and \(\Acal_{\mathrm{ALI}}\), are computationally efficient and easy to implement. In this section, we investigate their numerical performance and practical implications in a simulation study.

In our simulation, the target mean function is given by 
\[f^{[t]}(x) = x \cos(25x) + 4\abs{x-0.5}, \qquad(0\leq x\leq 1).\]
We also have two source samples and their mean functions are expressed in terms of difference functions:
\[ f^{[s,k]} = f^{[t]} - \delta^{[s,k]} ~~(k = 1,2), \qquad\text{where}~ \begin{cases}
    \delta^{[s,1]}(x) = -x^2, \\
    \delta^{[s,2]}(x) = e^x-1,
\end{cases} (0\leq x\leq 1).\]
It is worth noticing that the differences \(\delta^{[s,1]}\) and \(\delta^{[s,2]}\) are smoother than the target \(f^{[t]}\). This means that the differences are easier to estimate than the target, and transfer learning from the source samples can improve the estimation performance.

The target and source subjects in each sample have random curves generated as:
\begin{align*}
    X^{[t]}(x) &= f^{[t]}(x) + B_x^{[t]}, \\
    X^{[s,k]}(x) &= f^{[s,1]}(x) + B_x^{[s,k]} ~~(k=1,2),
\end{align*}
where \(\{B_x^{[t]}:x\geq 0\}\) and \(\{B_x^{[s,k]} : x\geq 0\}\) \((k=1,2)\) denote independent standard Brownian motions. It is well-known that Brownian motions have wriggly sample paths with smoothness strictly less than \(1/2\) almost surely. Traditional frameworks such as \citep{caiOptimalEstimationMean2011} do not allow for this type of model because the curves are not differentiable. However, the framework in this paper imposes no restrictions on the smoothness of the curves, making it more flexible and general even under conventional learning settings. Finally, each observation at a design point is assumed to be corrupted by additive standard Gaussian noise.

Let us first implement numerical experiments for a common design. The target and source samples consist of \(n_t\) and \(n_s\) subjects, respectively, and their curves are generated using the previously specified scheme. For each curve, equidistant observations at \(m_t\) and \(m_s\) design points, respectively, are taken on the domain \([0,1]\) with standard Gaussian noise. Following the estimation process in Theorem \ref{thm:upper bound under a common design}, the adaptive learning algorithm \(\Acal_{\mathrm{ALC}}\) for a common design is repeated \(r_{\mathrm{max}} = 20\) times and the final estimate \(\widehat{f}^{[t]}\) is obtained by taking an average of the results. We repeat this simulation 50 times, creating a boxplot of the \(\Lcal^2\)-integrated mean a squared loss to evaluate the risk for various combinations of \((n_t, m_t, n_s, m_s)\).

\begin{figure}[H]
    \centering
    \begin{subfigure}{.5\textwidth}
        \centering\captionsetup{justification=centering}
        \includegraphics[width=.95\linewidth]{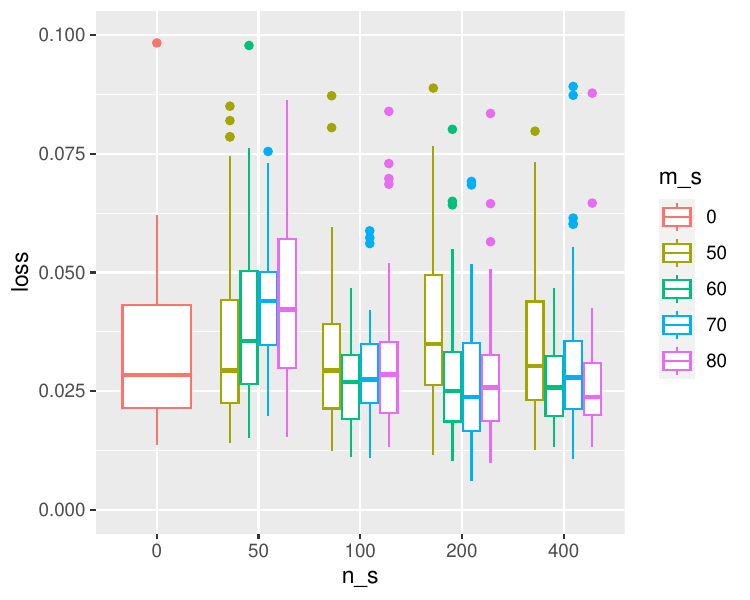}
        \caption{High-frequency regime\\(\(n_t = m_t = 50\))}
    \end{subfigure}%
    \begin{subfigure}{.5\textwidth}
        \centering\captionsetup{justification=centering}
        \includegraphics[width=.95\linewidth]{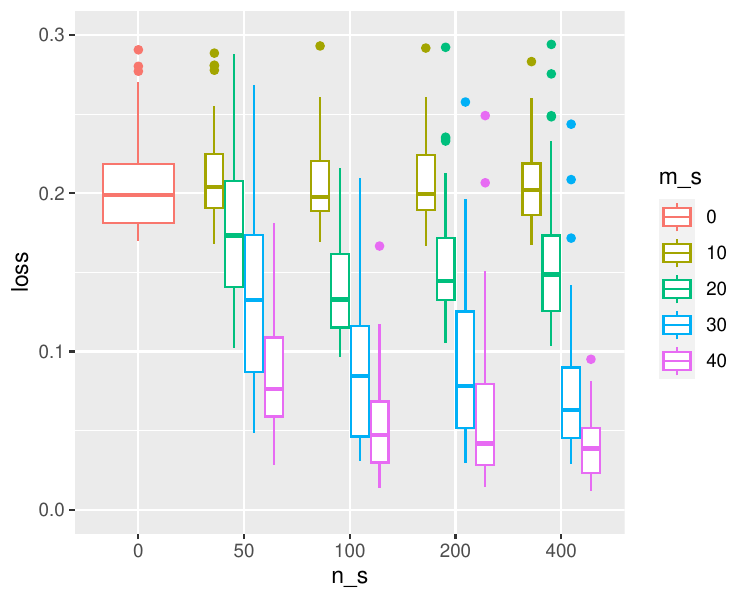}
        \caption{Low-frequency regime\\(\(n_t = 50,~ m_t = 10\))}
    \end{subfigure}%
    \caption{The results for numerical experiments under a common design}
    \label{fig:numerical experiments for a common design}
\end{figure}

Our theory indicates that a phase transition occurs at \(m_t \asymp n_t^{1/2\alpha_m}\), which determines the effectiveness of the source samples and transfer learning. To illustrate this phenomenon, we will compare two regimes: the high-frequency regime (\(n_t = m_t = 50\)) and the low-frequency regime (\(n_t = 50,~ m_t = 10\)). For the former regime, we shall consider all possible combinations between \(n_s \in \{50, 100, 200, 400\}\) and \(m_s \in \{50, 60, 70, 80\}\), whereas, for the latter regime, we consider \(m_s \in \{10, 20, 30, 40\}\) instead. Plus, we will consider the conventional learning setup \((n_s = m_s = 0)\) which serves as a baseline for any other setup. Figure \ref{fig:numerical experiments for a common design} shows the aggregated boxplots for all combinations in each regime.

The results from our theory are consistent with Figure \ref{fig:numerical experiments for a common design}. Transfer learning does not have a significant impact on the estimation performance in the high-frequency regime, while it is effective in the low-frequency regime. It is noteworthy that when \(m_s = 10\) in the low-frequency regime, the performance is similar to the conventional learning \((n_s = m_s = 0)\) regardless of the size \(n_s\) of the source subject. To enjoy benefits from the transfer learning, the number of design points for the source curves must exceed that of the target curves, i.e. \(m_t < m_s\).

\begin{figure}[H]
    \centering
    \begin{subfigure}{.5\textwidth}
        \centering\captionsetup{justification=centering}
        \includegraphics[width=.95\linewidth]{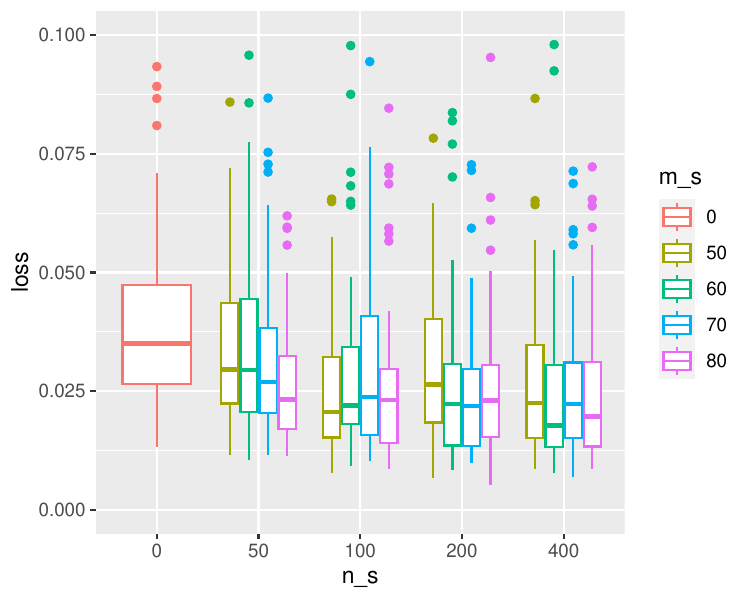}
        \caption{High-frequency regime\\(\(n_t = m_t = 50\))}
    \end{subfigure}%
    \begin{subfigure}{.5\textwidth}
        \centering\captionsetup{justification=centering}
        \includegraphics[width=.95\linewidth]{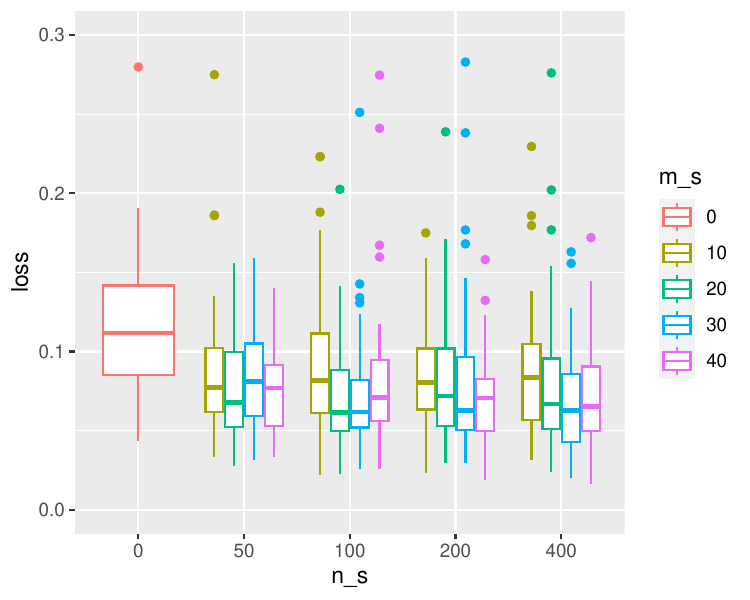}
        \caption{Low-frequency regime\\(\(n_t = 50,~ m_t = 10\))}
    \end{subfigure}%
        \caption{The results for numerical experiments under an independent design}
    \label{fig:numerical experiments for an independent design}
\end{figure}

We next concentrate on numerical experiments under an independent design. The target and source samples are generated using the same method as a common design, except that design points are independently and uniformly drawn from the domain \([0,1]\). We use the adaptive learning algorithm \(\Acal_{\mathrm{ALI}}\) instead of \(\Acal_{\mathrm{ALC}}\) for estimation while keeping all other simulation details the same. The resulting boxplots, similar to those in Figure \ref{fig:numerical experiments for a common design}, are presented in Figure \ref{fig:numerical experiments for an independent design}.

By comparing the high- and low-frequency regimes depicted in Figure \ref{fig:numerical experiments for an independent design}, we have come to the same conclusion as we did under a common design: transfer learning is only advantageous under the low-frequency regime. Additionally, we can draw interesting insights from comparing common and independent designs. For instance, under the high-frequency regime, both designs exhibit similar performance, which aligns with our theory as we always achieve the fastest parametric rate \(n_t^{-1}\) in this regime. Conversely, under the low-frequency regime, an independent design outperforms a common design when all other factors remain the same. To summarize, the results obtained from the numerical experiments are both explainable and consistent with our theory.

%% file: content/discussions.tex

\section{Discussion} \label{sec:discussion}

The present paper thoroughly investigated the transfer learning problem for functional mean estimation. It offers a comprehensive analysis of minimax convergence rates, which not only extends the previous results in conventional learning but also quantifies the advantages of transfer learning in terms of estimation performance. Our findings unveil a compelling phase transition phenomenon under both common and independent designs. Moreover, we developed data-driven algorithms that achieve optimal rates of convergence up to logarithmic factors across a broad spectrum of function spaces.

This paper serves as a critical stepping stone for further theoretical explorations in functional data analysis within the context of transfer learning. One promising avenue for future research is the investigation of transfer learning for the estimation of the covariance function, which plays a pivotal role in the analysis of functional and longitudinal data \citep{caiNonparametricCovarianceFunction2010}. In particular, the study of estimating the covariance function,
\[\Sigma^{[t]}(u,v) := \Ebb \Bigl[\bigl(X^{[t]}(u) - \Ebb X^{[t]}(u)\bigr) \bigl(X^{[t]}(v) - \Ebb X^{[t]}(v)\bigr)\Bigr], \qquad(u,v \in [0,1]), \]
should be of substantial interest, given the availability of both target and source samples, as presented in Equations \eqref{eq:modeling for target sample} and \eqref{eq:modeling for source samples}. Establishing the minimax rates of convergence for this problem promises to yield valuable insights.

Another important problem is transfer learning for functional linear regression. In the conventional settings, \citet{yuanReproducingKernelHilbert2010} established the minimax rate of convergence for estimating the regression function, while both \citet{caiPredictionFunctionalLinear2006} and \citet{caiMinimaxAdaptivePrediction2012} investigated the optimal prediction problem. It is of significant interest to explore transfer learning for linear regression and investigate the conditions under which source samples can improve the estimation performance.

While our present work offers valuable insights, it is essential to acknowledge its limitations and explore avenues for further improvement. One notable constraint is our assumption that both the target and source samples share the same sampling design. However, in practice, the designs can differ from each other. For example,  the target sample can be generated from a common design while the source samples are obtained through an independent design. Not only that, some source samples are drawn from a common design, while others come from an independent design. Even within source samples sharing the same design, there could be variations in the common design points or the common marginal distributions that generate design points. 

An analogous assumption is imposed on the smoothness of mean functions as shown in Equation \eqref{eq:modeling assumption for mean and difference functions}, but this could be generalized to cover more diverse datasets. For example, we could assume some source mean functions are smoother while others are rougher than the target: for given \(K_1, K_2 \in \Zbb^+\) with \(K_1 + K_2 = K\),
\[\begin{aligned}
    &f^{[s,1]}, \ldots, f^{[s,K_1]} \in \Hcal_{\alpha_{m_1}}(L_{m_1}, M_{m_1}),\\
    &f^{[s,K_1+1]}, \ldots, f^{[s,K]} \in \Hcal_{\alpha_{m_2}}(L_{m_2}, M_{m_2}),
\end{aligned} \qquad\text{for some}~ \alpha_{m_1}, \alpha_{m_2} > 0. \]

Another aspect to consider is the assumption that all source samples share the same number of subjects and design points, as presented in Equation \eqref{eq:modeling assumption for the number of subjects and design points}. However, this assumption can be relaxed to accommodate different numbers of subjects or design points across the source samples. While this relaxation may lead to a more complex minimax rate of convergence, it would better align with the diversity observed in real-world datasets.

Finally, in the independent design setting, we assumed that the marginal distribution of the source samples only trivially differs from that of the target sample. In practice, the marginal density ratio between the target and source samples may not be uniformly bounded away from zero and infinity. This departure from the uniformity of density ratio can significantly impact the minimax rate of convergence. For example, if the marginal distributions of the target and source samples are completely singular, the benefit of transfer learning may be negligible. Incorporating this effect into the analysis could be achieved by introducing a new measure to quantify the singularity between marginal distributions. Similar to the transfer exponent introduced in \citet{kpotufeMarginalSingularityBenefits2018}, such a measure would be a valuable addition to assess the influence of marginal singularity under the transfer learning setup.

%% file: content/proofs.tex

\section{Proofs} \label{sec:proof}

This section contains the proofs of our main results for common design, especially under transfer learning setup: Theorem \ref{thm:upper bound under a common design} and Theorem \ref{thm:lower bound under a common design}. We will defer the proofs of the remaining results for common design, as well as the results for independent design and every other technical results, to the Supplementary Material \citep{caiSupplementTransferLearning2024}.

\subsection{Proof of Theorem \ref{thm:upper bound under a common design}}

Thanks to Theorem \ref{thm:minimax risk under conventional setup and common design} and the nature of our estimator \(\widehat{f}^{[t]}\), it suffices to argue the following rate of the transfer learning estimator \(\widehat{f}_{\mathrm{TL}}^{[t]}\),
\[ \sup_{\Pbb\in\Pcal} \Ebb\norm{\widehat{f}_{\mathrm{TL}}^{[t]} - f^{[t]}}^2_{\Lcal^2} \lesssim L_\delta^2 m_t^{-2\alpha_\delta} + \frac{\log^2 n_t n_s}{n_t} + L_m^2 m_s^{-2\alpha_m} + \frac{\log^2 (K n_s)}{K n_s}. \]
Using the notations in the transfer learning algorithm \(\Acal_{\mathrm{TL}}\), we may easily check
\[ \sup_{\Pbb\in\Pcal} \Ebb\norm{\widehat{f}_{\mathrm{TL}}^{[t]} - f^{[t]}}^2_{\Lcal^2} \leq 2 \sup_{\Pbb\in\Pcal} \Ebb\norm{\widehat{f}^{[s]} - f^{[s]}}^2_{\Lcal^2} + 2 \sup_{\Pbb\in\Pcal} \Ebb\norm{\widehat{\delta}^{[s]} - \delta^{[s]}}^2_{\Lcal^2} \]
We shall present the rates of convergence for those estimators, \(\widehat{f}^{[s]}\) and \(\widehat{\delta}^{[s]}\), in the following propositions. The proofs are given in the Supplementary material \citep{caiSupplementTransferLearning2024}.

\begin{proposition}
    If \(m_s b_s \geq 2C_s(d_s+1)\) and \(\omega(\alpha_m) \leq d_s \leq O(1)\) are satisfied,
    \[ \sup_{\Pbb\in\Pcal} \Ebb\norm{\widehat{f}^{[s]} - f^{[s]}}^2_{\Lcal^2} \lesssim L_m^2 b_s^{2\alpha_m} + \frac{\log^2 (K n_s)}{K n_s} \]
\end{proposition}

\begin{proposition}
    If \(m_t b_\delta \geq 2C_t(d_\delta+1)\) and \(\omega(\alpha_\delta) \leq d_\delta \leq O(1)\) are fulfilled,
    \[ \sup_{\Pbb\in\Pcal} \Ebb\norm{\widehat{\delta}^{[s]} - \delta^{[s]}}^2_{\Lcal^2} \lesssim L_\delta^2 b_\delta^{2\alpha_\delta} + \frac{\log^2 n_t n_s}{n_t}. \]
\end{proposition}

It is therefore optimal to choose the bandwidths \(b_s = \Theta(\inverse{m_s})\) and \(b_\delta = \Theta(\inverse{m_t})\) such that they satisfy \(m_s b_s \geq 2C_s(d_s+1)\) and \(m_t b_\delta \geq 2C_t(d_\delta+1)\), and to select any constant degrees \(d_s \geq \omega(\alpha_m)\) and \(d_\delta \geq \omega(\alpha_\delta)\). This immediately concludes the proof.
\qed

\subsection{Proof of Theorem \ref{thm:lower bound under a common design}}

It suffices to show that for any estimator \(\widehat{f}^{[t]}\),
\[ \begin{aligned}
    \sup_{\Pbb\in\Pcal} \Ebb\norm{\widehat{f}^{[t]} - f^{[t]}}^2_{\Lcal^2} &\gtrsim \bigl(L_\delta^2 m_t^{-2\alpha_\delta} \binmin L_m^2 m_t^{-2\alpha_m}\bigr) + \frac{1}{n_t} \\
    &~~~~+ \left(L_m^2 m_t^{-2\alpha_m} + \frac{1}{n_t} \right) \binmin \left(L_m^2 m_s^{-2\alpha_m} + \frac{1}{K n_s}\right).
\end{aligned} \]

Consider a sub-model \(\Pcal_1 \subset \Pcal\) where the source samples \(\Dcal^{[s,1]}, \ldots, \Dcal^{[s,K]}\) are completely futile in the sense that \(f^{[s,1]} = \cdots = f^{[s,K]} = 0\). Under this sub-model \(\Pcal_1\), target sample \(\Dcal^{[t]}\) forms a sufficient statistic of the target mean function \(f^{[t]}\). On top of that, the target mean function should be further restricted by \(f^{[t]} \in \Hcal_{\alpha_m}(L_m, M_m) \cap \Hcal_{\alpha_\delta}(L_\delta, M_\delta)\). As a result, Theorem \ref{thm:minimax risk under conventional setup and common design} immediately implies that:
\[ \sup_{\Pbb\in\Pcal} \Ebb\norm{\widehat{f}^{[t]} - f^{[t]}}^2_{\Lcal^2} \geq \sup_{\Pbb\in\Pcal_1} \Ebb\normbig{\Ebb(\widehat{f}^{[t]} \| \Dcal^{[t]}) - f^{[t]}}^2_{\Lcal^2} \gtrsim \bigl(L_\delta^2 m_t^{-2\alpha_\delta} \binmin L_m^2 m_t^{-2\alpha_m}\bigr) + \frac{1}{n_t}. \]

We are now enough to argue the other part of the lower bound:
\[ \sup_{\Pbb\in\Pcal} \Ebb\norm{\widehat{f}^{[t]} - f^{[t]}}^2_{\Lcal^2} \gtrsim \left(L_m^2 m_t^{-2\alpha_m} + \frac{1}{n_t} \right) \binmin \left(L_m^2 m_s^{-2\alpha_m} + \frac{1}{K n_s}\right). \]
This is equivalent to insisting on the subsequent three propositions:

\begin{proposition} \label{prop:lower bound under a common design case 1}
    For any estimator \(\widehat{f}^{[t]}\), we have
    \[\sup_{\Pbb\in\Pcal} \Ebb\norm{\widehat{f}^{[t]} - f^{[t]}}^2_{\Lcal^2} \gtrsim \frac{1}{n_t} \binmin \frac{1}{K n_s}.\]
\end{proposition}

\begin{proposition} \label{prop:lower bound under a common design case 2}
    For any estimator \(\widehat{f}^{[t]}\), we have
    \[ \sup_{\Pbb\in\Pcal} \Ebb\norm{\widehat{f}^{[t]} - f^{[t]}}^2_{\Lcal^2} \gtrsim L_m^2 \bigl(m_t^{-2\alpha_m} \binmin  m_s^{-2\alpha_m}\bigr). \]
\end{proposition}

\begin{proposition} \label{prop:lower bound under a common design case 3}
    For any estimator \(\widehat{f}^{[t]}\), we have
    \begin{align*}
        \sup_{\Pbb\in\Pcal} \Ebb\norm{\widehat{f}^{[t]} - f^{[t]}}^2_{\Lcal^2} \gtrsim L_m^2 m_t^{-2\alpha_m} \binmin \frac{1}{K n_s}.
    \end{align*}
\end{proposition}

\begin{proposition} \label{prop:lower bound under a common design case 4}
    For any estimator \(\widehat{f}^{[t]}\), we have
    \begin{align*}
        \sup_{\Pbb\in\Pcal} \Ebb\norm{\widehat{f}^{[t]} - f^{[t]}}^2_{\Lcal^2} \gtrsim L_m^2 m_s^{-2\alpha_m} \binmin \frac{1}{n_t}.
    \end{align*}
\end{proposition}

We construct another sub-model \(\Pcal_2 \subset \Pcal\) as follows. The source samples \(\Dcal^{[s,1]}, \ldots, \Dcal^{[s,K]}\) are the most advantageous in the sense that \(f^{[t]} = f^{[s,1]} = \cdots = f^{[s,K]}\). In other words, every sample \(\Dcal^{[t]}, \Dcal^{[s,1]}, \ldots, \Dcal^{[s,K]}\) is assumed to share the same data-generating process whose further specifications would be different depending on each case. Let us take any \(C^{\infty}\) function \(\phi:\Rbb\to\Rbb\) such that \(\phi^{(q)}(x)=0\) for any \(q=0,1,2,\ldots\) and \(x\not\in (0,1)\). We will prove Proposition \ref{prop:lower bound under a common design case 1} and \ref{prop:lower bound under a common design case 2}. The proof of Proposition \ref{prop:lower bound under a common design case 3} and \ref{prop:lower bound under a common design case 4} will be presented separately in the Supplementary material \citep{caiSupplementTransferLearning2024}. \qed

\subsubsection*{Proof of Proposition \ref{prop:lower bound under a common design case 1}}

Under the sub-model \(\Pcal_2 \subset \Pcal\), let us assume the target random function \(X^{[t]}_1\) as well as the source ones \(X^{[s,1]}_{1}, \ldots, X^{[s,K]}_{1}\) are unknown constant functions. Since their mean functions should coincide, the problem is equivalent to estimating the scalar mean from \(n_t + K n_s\) number of independent and identically distributed observations under mean squared error. It is thus immediate to get the lower bound of the parametric rate:
\[ \sup_{\Pbb\in\Pcal} \Ebb\norm{\widehat{f}^{[t]} - f^{[t]}}^2_{\Lcal^2} \geq \sup_{\Pbb\in\Pcal_2} \Ebb\norm{\widehat{f}^{[t]} - f^{[t]}}^2_{\Lcal^2} \gtrsim \frac{1}{n_t + K n_s} \gtrsim \frac{1}{n_t} \binmin \frac{1}{K n_s}. \]
\qed

\subsubsection*{Proof of Proposition \ref{prop:lower bound under a common design case 2}}

We begin by letting a quantity \(m:= m_t + m_s\). For each vector \(v = (v_1, \ldots, v_{2m}) \in \{0,1\}^{2m}\) on hypercube, we define a function \(g_v : [0,1] \to \Rbb\) by
\[ g_v(x) := \sum_{r=1}^{2m} C L_m v_r (2m)^{-\alpha_m} \phi\bigl(2m x - (r-1)\bigr), \]
where \(C > 0\) is a small enough generic constant such that \(g_v \in \Hcal_{\alpha_m}(L_m, M_m)\). The sub-model \(\Pcal_2 \subset \Pcal\) consists of any probability measure \(\Pbb_v\) \((v \in \{0,1\}^{2m})\) whose common mean function \(f^{[t]} = f^{[s,1]}\) is given as \(g_v\). This sub-model \(\Pcal_2\) is Hamming separated in that
\[\norm{g_{v} - g_{w}}_{\Lcal^2}^2 \geq C^2 L_m^2 (2m)^{-2\alpha_m - 1} \norm{\phi}_{\Lcal^2}^2 \sum_{r=1}^{2m} \mathbb1 (v_r \not= w_r), \quad\text{for any}~ v, w \in \{0,1\}^{2m}. \]
For each \(r = 1,\ldots,2m\), let \(\overline{\Pbb}_{0,r}\) and \(\overline{\Pbb}_{1,r}\) denote the mixture distributions of \(\{\Pbb_v : v_r = 0\}\) and \(\{\Pbb_v : v_r = 1\}\), respectively. It follows from Assouad's lemma \citep{assouadDensiteDimension1983} that
\begin{align*}
    \sup_{\Pbb\in\Pcal} \Ebb\norm{\widehat{f}^{[t]} - f^{[t]}}^2_{\Lcal^2} &\geq \sup_{\Pbb\in\Pcal_2} \Ebb\norm{\widehat{f}^{[t]} - f^{[t]}}^2_{\Lcal^2} \numberthis\label{eq:application of assouad's lemma under a common design case 1} \\
    &\geq \frac{1}{2} C^2 L_m^2 (2m)^{-2\alpha_m - 1} \norm{\phi}_{\Lcal^2}^2 \sum_{r=1}^{2m} \Big( 1 - \norm{\overline{\Pbb}_{0,r} - \overline{\Pbb}_{1,r}}_{\mathrm{TV}} \Big),
\end{align*}
where \(\norm{\,\cdot\,}_{\mathrm{TV}}\) denotes the total variation distance. Notice that two mixtures \(\overline{\Pbb}_{0,r}\) and \(\overline{\Pbb}_{1,r}\) follow exactly the same distribution when the interval \([(r-1)/2m, r/2m)\) contains nothing among \(\{T_{j}^{[t]}:j=1,\ldots,m_t\} \cup \{T_{j}^{[s,1]}:j=1,\ldots,m_s\}\), the common design points for target and source domains. Since these design points are at most \(m = m_t + m_s\), the number of such an \(r \in \{1,\ldots,2m\}\) must be at least \(m\) by Pigeonhole principle. In other words, we have
\[ \sum_{r=1}^{2m} \Big( 1 - \norm{\overline{\Pbb}_{0,r} - \overline{\Pbb}_{1,r}}_{\mathrm{TV}} \Big) \geq m. \numberthis\label{eq:bound for total variation distance under a common design case 1} \]
The desired result is now immediate by combining Equation \eqref{eq:application of assouad's lemma under a common design case 1} and \eqref{eq:bound for total variation distance under a common design case 1}:
\[ \pushQED{\qed} \sup_{\Pbb\in\Pcal} \Ebb\norm{\widehat{f}^{[t]} - f^{[t]}}^2_{\Lcal^2} \gtrsim L_m^2 m^{-2\alpha_m} \asymp L_m^2 \bigl(m_t^{-2\alpha_m} \binmin m_s^{-2\alpha_m}\bigr). \qedhere\popQED \]

%% file: content/epilogue.tex

\begin{acks}[Acknowledgments]
    The authors extend their gratitude to the Editor, the Associate Editor and the anonymous referees for their constructive and insightful comments and suggestions, which have significantly enhanced the quality of this paper.
\end{acks}

\begin{funding}
The research of Tony Cai was supported in part by NSF Grant DMS-2015259 and NIH grant R01-GM129781. 
\end{funding}

\begin{supplement}
    \stitle{Supplement to ``Transfer Learning for Functional Mean Estimation: Phase Transition and Adaptive Algorithms.''} 
    \slink[doi]{0000000} \sdatatype{.pdf}
    \sdescription{This supplementary material provides the complete proofs of the main theorems and technical results introduced in the paper, ``Transfer Learning for Functional Mean Estimation: Phase Transition and Adaptive Algorithms.'' The structure of this material is as follows. Supplementary Material A handles the remaining results in common design: Theorem 2.1 and 2.4. Supplementary Material B covers the main results in independent design: Theorems 3.1, 3.2, 3.3, and 3.4. Supplementary Material C includes the proofs of technical results relevant to common design such as Proposition 6.1, 6.2, 6.5 and 6.6 as well as Lemma A.1 and C.1. Finally, Supplementary Material D contains the proofs of Lemma B.1--B.5, which are technical lemmas for the main results in independent design.}
\end{supplement}

%% file: Transfer-Learning-for-Mean-Function-Estimation.bbl
\begin{thebibliography}{42}

\bibitem[\protect\citeauthoryear{Assouad}{1983}]{assouadDensiteDimension1983}
\begin{barticle}[author]
\bauthor{\bsnm{Assouad},~\bfnm{Patrick}\binits{P.}}
(\byear{1983}).
\btitle{Densit{\'e} et Dimension}.
\bjournal{Annales de l'Institut Fourier}
\bvolume{33}
\bpages{233--282}.
\bdoi{10.5802/aif.938}
\end{barticle}
\endbibitem

\bibitem[\protect\citeauthoryear{Breiman}{1996}]{breimanBaggingPredictors1996}
\begin{barticle}[author]
\bauthor{\bsnm{Breiman},~\bfnm{Leo}\binits{L.}}
(\byear{1996}).
\btitle{Bagging Predictors}.
\bjournal{Machine Learning}
\bvolume{24}
\bpages{123--140}.
\bdoi{10.1007/BF00058655}
\end{barticle}
\endbibitem

\bibitem[\protect\citeauthoryear{Cai, Cai and
  Li}{2022}]{caiTransferLearningContextual2022}
\begin{barticle}[author]
\bauthor{\bsnm{Cai},~\bfnm{Changxiao}\binits{C.}},
  \bauthor{\bsnm{Cai},~\bfnm{T.~Tony}\binits{T.~T.}} \AND
  \bauthor{\bsnm{Li},~\bfnm{Hongzhe}\binits{H.}}
(\byear{2022}).
\btitle{Transfer {{Learning}} for {{Contextual Multi-armed Bandits}}}.
\bjournal{Technical report}.
\bdoi{10.48550/arXiv.2211.12612}
\end{barticle}
\endbibitem

\bibitem[\protect\citeauthoryear{Cai and
  Hall}{2006}]{caiPredictionFunctionalLinear2006}
\begin{barticle}[author]
\bauthor{\bsnm{Cai},~\bfnm{T.~Tony}\binits{T.~T.}} \AND
  \bauthor{\bsnm{Hall},~\bfnm{Peter}\binits{P.}}
(\byear{2006}).
\btitle{Prediction in Functional Linear Regression}.
\bjournal{The Annals of Statistics}
\bvolume{34}
\bpages{2159--2179}.
\bdoi{10.1214/009053606000000830}
\end{barticle}
\endbibitem

\bibitem[\protect\citeauthoryear{Cai, Kim and
  Pu}{2024}]{caiSupplementTransferLearning2024}
\begin{barticle}[author]
\bauthor{\bsnm{Cai},~\bfnm{T.~Tony}\binits{T.~T.}},
  \bauthor{\bsnm{Kim},~\bfnm{Dongwoo}\binits{D.}} \AND
  \bauthor{\bsnm{Pu},~\bfnm{Hongming}\binits{H.}}
(\byear{2024}).
\btitle{Supplement to ``{{Transfer Learning}} for {{Functional Mean
  Estimation}}: {{Phase Transition}} and {{Adaptive Algorithms}}''}.
\end{barticle}
\endbibitem

\bibitem[\protect\citeauthoryear{Cai and
  Pu}{2022}]{caiTransferLearningNonparametric2022}
\begin{barticle}[author]
\bauthor{\bsnm{Cai},~\bfnm{T.~Tony}\binits{T.~T.}} \AND
  \bauthor{\bsnm{Pu},~\bfnm{Hongming}\binits{H.}}
(\byear{2022}).
\btitle{Transfer Learning for Nonparametric Regression: Non-Asymptotic Minimax
  Analysis and Adaptive Procedure}.
\bjournal{Technical report}.
\end{barticle}
\endbibitem

\bibitem[\protect\citeauthoryear{Cai and
  Wei}{2021}]{caiTransferLearningNonparametric2021}
\begin{barticle}[author]
\bauthor{\bsnm{Cai},~\bfnm{T.~Tony}\binits{T.~T.}} \AND
  \bauthor{\bsnm{Wei},~\bfnm{Hongji}\binits{H.}}
(\byear{2021}).
\btitle{Transfer Learning for Nonparametric Classification: {{Minimax}} Rate
  and Adaptive Classifier}.
\bjournal{Annals of Statistics}
\bvolume{49}
\bpages{100--128}.
\bdoi{10.1214/20-AOS1949}
\end{barticle}
\endbibitem

\bibitem[\protect\citeauthoryear{Cai and
  Yuan}{2010}]{caiNonparametricCovarianceFunction2010}
\begin{barticle}[author]
\bauthor{\bsnm{Cai},~\bfnm{T.~Tony}\binits{T.~T.}} \AND
  \bauthor{\bsnm{Yuan},~\bfnm{Ming}\binits{M.}}
(\byear{2010}).
\btitle{Nonparametric Covariance Function Estimation for Functional and
  Longitudinal Data}.
\bjournal{Technical report}.
\end{barticle}
\endbibitem

\bibitem[\protect\citeauthoryear{Cai and
  Yuan}{2011}]{caiOptimalEstimationMean2011}
\begin{barticle}[author]
\bauthor{\bsnm{Cai},~\bfnm{T.~Tony}\binits{T.~T.}} \AND
  \bauthor{\bsnm{Yuan},~\bfnm{Ming}\binits{M.}}
(\byear{2011}).
\btitle{Optimal Estimation of the Mean Function Based on Discretely Sampled
  Functional Data: {{Phase}} Transition}.
\bjournal{The Annals of Statistics}
\bvolume{39}
\bpages{2330--2355}.
\bdoi{10.1214/11-AOS898}
\end{barticle}
\endbibitem

\bibitem[\protect\citeauthoryear{Cai and
  Yuan}{2012}]{caiMinimaxAdaptivePrediction2012}
\begin{barticle}[author]
\bauthor{\bsnm{Cai},~\bfnm{T.~Tony}\binits{T.~T.}} \AND
  \bauthor{\bsnm{Yuan},~\bfnm{Ming}\binits{M.}}
(\byear{2012}).
\btitle{Minimax and {{Adaptive Prediction}} for {{Functional Linear
  Regression}}}.
\bjournal{Journal of the American Statistical Association}
\bvolume{107}
\bpages{1201--1216}.
\bdoi{10.1080/01621459.2012.716337}
\end{barticle}
\endbibitem

\bibitem[\protect\citeauthoryear{Choi
  et~al.}{2017}]{choiTransferLearningMusic2017}
\begin{barticle}[author]
\bauthor{\bsnm{Choi},~\bfnm{Keunwoo}\binits{K.}},
  \bauthor{\bsnm{Fazekas},~\bfnm{Gy{\"o}rgy}\binits{G.}},
  \bauthor{\bsnm{Sandler},~\bfnm{Mark}\binits{M.}} \AND
  \bauthor{\bsnm{Cho},~\bfnm{Kyunghyun}\binits{K.}}
(\byear{2017}).
\btitle{Transfer Learning for Music Classification and Regression Tasks}.
\bjournal{Technical report}.
\bdoi{10.48550/arXiv.1703.09179}
\end{barticle}
\endbibitem

\bibitem[\protect\citeauthoryear{Degras}{2017}]{degrasSimultaneousConfidenceBands2017}
\begin{barticle}[author]
\bauthor{\bsnm{Degras},~\bfnm{David}\binits{D.}}
(\byear{2017}).
\btitle{Simultaneous Confidence Bands for the Mean of Functional Data}.
\bjournal{WIREs Computational Statistics}
\bvolume{9}
\bpages{e1397}.
\bdoi{10.1002/wics.1397}
\end{barticle}
\endbibitem

\bibitem[\protect\citeauthoryear{Gong
  et~al.}{2012}]{gongGeodesicFlowKernel2012}
\begin{binproceedings}[author]
\bauthor{\bsnm{Gong},~\bfnm{Boqing}\binits{B.}},
  \bauthor{\bsnm{Shi},~\bfnm{Yuan}\binits{Y.}},
  \bauthor{\bsnm{Sha},~\bfnm{Fei}\binits{F.}} \AND
  \bauthor{\bsnm{Grauman},~\bfnm{Kristen}\binits{K.}}
(\byear{2012}).
\btitle{Geodesic Flow Kernel for Unsupervised Domain Adaptation}.
In \bbooktitle{2012 {{IEEE Conference}} on {{Computer Vision}} and {{Pattern
  Recognition}}}
\bpages{2066--2073}.
\bdoi{10.1109/CVPR.2012.6247911}
\end{binproceedings}
\endbibitem

\bibitem[\protect\citeauthoryear{Huang
  et~al.}{2013}]{huangCrosslanguageKnowledgeTransfer2013}
\begin{binproceedings}[author]
\bauthor{\bsnm{Huang},~\bfnm{Jui-Ting}\binits{J.-T.}},
  \bauthor{\bsnm{Li},~\bfnm{Jinyu}\binits{J.}},
  \bauthor{\bsnm{Yu},~\bfnm{Dong}\binits{D.}},
  \bauthor{\bsnm{Deng},~\bfnm{Li}\binits{L.}} \AND
  \bauthor{\bsnm{Gong},~\bfnm{Yifan}\binits{Y.}}
(\byear{2013}).
\btitle{Cross-Language Knowledge Transfer Using Multilingual Deep Neural
  Network with Shared Hidden Layers}.
In \bbooktitle{2013 {{IEEE International Conference}} on {{Acoustics}},
  {{Speech}} and {{Signal Processing}}}
\bpages{7304--7308}.
\bdoi{10.1109/ICASSP.2013.6639081}
\end{binproceedings}
\endbibitem

\bibitem[\protect\citeauthoryear{James and
  Menzies}{2021}]{jamesTrendsCOVID19Prevalence2021}
\begin{barticle}[author]
\bauthor{\bsnm{James},~\bfnm{Nick}\binits{N.}} \AND
  \bauthor{\bsnm{Menzies},~\bfnm{Max}\binits{M.}}
(\byear{2021}).
\btitle{Trends in {{COVID-19}} Prevalence and Mortality: {{A}} Year in Review}.
\bjournal{Physica D: Nonlinear Phenomena}
\bvolume{425}
\bpages{132968}.
\bdoi{10.1016/j.physd.2021.132968}
\end{barticle}
\endbibitem

\bibitem[\protect\citeauthoryear{Jiang, Aston and
  Wang}{2009}]{jiangSmoothingDynamicPositron2009}
\begin{barticle}[author]
\bauthor{\bsnm{Jiang},~\bfnm{Ci-Ren}\binits{C.-R.}},
  \bauthor{\bsnm{Aston},~\bfnm{John A.~D.}\binits{J.~A.~D.}} \AND
  \bauthor{\bsnm{Wang},~\bfnm{Jane-Ling}\binits{J.-L.}}
(\byear{2009}).
\btitle{Smoothing Dynamic Positron Emission Tomography Time Courses Using
  Functional Principal Components}.
\bjournal{NeuroImage}
\bvolume{47}
\bpages{184--193}.
\bdoi{10.1016/j.neuroimage.2009.03.051}
\end{barticle}
\endbibitem

\bibitem[\protect\citeauthoryear{Kalogridis and
  Van~Aelst}{2022}]{kalogridisRobustOptimalEstimation2022}
\begin{barticle}[author]
\bauthor{\bsnm{Kalogridis},~\bfnm{Ioannis}\binits{I.}} \AND
  \bauthor{\bsnm{Van~Aelst},~\bfnm{Stefan}\binits{S.}}
(\byear{2022}).
\btitle{Robust Optimal Estimation of Location from Discretely Sampled
  Functional Data}.
\bjournal{Scandinavian Journal of Statistics}.
\bdoi{10.1111/sjos.12586}
\end{barticle}
\endbibitem

\bibitem[\protect\citeauthoryear{Kozloff
  et~al.}{2020}]{kozloffCOVID19GlobalPandemic2020}
\begin{barticle}[author]
\bauthor{\bsnm{Kozloff},~\bfnm{Nicole}\binits{N.}},
  \bauthor{\bsnm{Mulsant},~\bfnm{Benoit~H}\binits{B.~H.}},
  \bauthor{\bsnm{Stergiopoulos},~\bfnm{Vicky}\binits{V.}} \AND
  \bauthor{\bsnm{Voineskos},~\bfnm{Aristotle~N}\binits{A.~N.}}
(\byear{2020}).
\btitle{The {{COVID-19 Global Pandemic}}: {{Implications}} for {{People With
  Schizophrenia}} and {{Related Disorders}}}.
\bjournal{Schizophrenia Bulletin}
\bvolume{46}
\bpages{752--757}.
\bdoi{10.1093/schbul/sbaa051}
\end{barticle}
\endbibitem

\bibitem[\protect\citeauthoryear{Kpotufe and
  Martinet}{2018}]{kpotufeMarginalSingularityBenefits2018}
\begin{binproceedings}[author]
\bauthor{\bsnm{Kpotufe},~\bfnm{Samory}\binits{S.}} \AND
  \bauthor{\bsnm{Martinet},~\bfnm{Guillaume}\binits{G.}}
(\byear{2018}).
\btitle{Marginal {{Singularity}}, and the {{Benefits}} of {{Labels}} in
  {{Covariate-Shift}}}.
In \bbooktitle{Proceedings of the 31st {{Conference On Learning Theory}}}
\bpages{1882--1886}.
\bpublisher{{PMLR}}.
\end{binproceedings}
\endbibitem

\bibitem[\protect\citeauthoryear{Leng and
  M{\"u}ller}{2006}]{lengClassificationUsingFunctional2006}
\begin{barticle}[author]
\bauthor{\bsnm{Leng},~\bfnm{Xiaoyan}\binits{X.}} \AND
  \bauthor{\bsnm{M{\"u}ller},~\bfnm{Hans-Georg}\binits{H.-G.}}
(\byear{2006}).
\btitle{Classification Using Functional Data Analysis for Temporal Gene
  Expression Data}.
\bjournal{Bioinformatics}
\bvolume{22}
\bpages{68--76}.
\bdoi{10.1093/bioinformatics/bti742}
\end{barticle}
\endbibitem

\bibitem[\protect\citeauthoryear{Li, Cai and
  Li}{2022a}]{liTransferLearningHighDimensional2022}
\begin{barticle}[author]
\bauthor{\bsnm{Li},~\bfnm{Sai}\binits{S.}},
  \bauthor{\bsnm{Cai},~\bfnm{T.~Tony}\binits{T.~T.}} \AND
  \bauthor{\bsnm{Li},~\bfnm{Hongzhe}\binits{H.}}
(\byear{2022}a).
\btitle{Transfer {{Learning}} for {{High-Dimensional Linear Regression}}:
  {{Prediction}}, {{Estimation}} and {{Minimax Optimality}}}.
\bjournal{Journal of the Royal Statistical Society Series B: Statistical
  Methodology}
\bvolume{84}
\bpages{149--173}.
\bdoi{10.1111/rssb.12479}
\end{barticle}
\endbibitem

\bibitem[\protect\citeauthoryear{Li, Cai and
  Li}{2022b}]{liTransferLearningLargeScale2022}
\begin{barticle}[author]
\bauthor{\bsnm{Li},~\bfnm{Sai}\binits{S.}},
  \bauthor{\bsnm{Cai},~\bfnm{T.~Tony}\binits{T.~T.}} \AND
  \bauthor{\bsnm{Li},~\bfnm{Hongzhe}\binits{H.}}
(\byear{2022}b).
\btitle{Transfer {{Learning}} in {{Large-Scale Gaussian Graphical Models}} with
  {{False Discovery Rate Control}}}.
\bjournal{Journal of the American Statistical Association}
\bvolume{0}
\bpages{1--13}.
\bdoi{10.1080/01621459.2022.2044333}
\end{barticle}
\endbibitem

\bibitem[\protect\citeauthoryear{Li
  et~al.}{2023}]{liEstimationInferenceHighDimensional2023}
\begin{barticle}[author]
\bauthor{\bsnm{Li},~\bfnm{Sai}\binits{S.}},
  \bauthor{\bsnm{Zhang},~\bfnm{Linjun}\binits{L.}},
  \bauthor{\bsnm{Cai},~\bfnm{T.~Tony}\binits{T.~T.}} \AND
  \bauthor{\bsnm{Li},~\bfnm{Hongzhe}\binits{H.}}
(\byear{2023}).
\btitle{Estimation and {{Inference}} for {{High-Dimensional Generalized Linear
  Models}} with {{Knowledge Transfer}}}.
\bjournal{Journal of the American Statistical Association}
\bvolume{0}
\bpages{1--12}.
\bdoi{10.1080/01621459.2023.2184373}
\end{barticle}
\endbibitem

\bibitem[\protect\citeauthoryear{Mant{\'e}, Durbec and
  Dauvin}{2005}]{manteFunctionalDataanalyticApproach2005}
\begin{barticle}[author]
\bauthor{\bsnm{Mant{\'e}},~\bfnm{C.}\binits{C.}},
  \bauthor{\bsnm{Durbec},~\bfnm{J.~p.}\binits{J.~p.}} \AND
  \bauthor{\bsnm{Dauvin},~\bfnm{J.~c.}\binits{J.~c.}}
(\byear{2005}).
\btitle{A Functional Data-Analytic Approach to the Classification of Species
  According to Their Spatial Dispersion. {{Application}} to a Marine
  Macrobenthic Community from the {{Bay}} of {{Morlaix}} ({{Western English
  Channel}})}.
\bjournal{Journal of Applied Statistics}
\bvolume{32}
\bpages{831--840}.
\bdoi{10.1080/02664760500080124}
\end{barticle}
\endbibitem

\bibitem[\protect\citeauthoryear{Page
  et~al.}{2006}]{pageNormalizingTemporalPatterns2006}
\begin{barticle}[author]
\bauthor{\bsnm{Page},~\bfnm{A.}\binits{A.}},
  \bauthor{\bsnm{Ayala},~\bfnm{G.}\binits{G.}},
  \bauthor{\bsnm{Le{\'o}n},~\bfnm{M.~T.}\binits{M.~T.}},
  \bauthor{\bsnm{Peydro},~\bfnm{M.~F.}\binits{M.~F.}} \AND
  \bauthor{\bsnm{Prat},~\bfnm{J.~M.}\binits{J.~M.}}
(\byear{2006}).
\btitle{Normalizing Temporal Patterns to Analyze Sit-to-Stand Movements by
  Using Registration of Functional Data}.
\bjournal{Journal of Biomechanics}
\bvolume{39}
\bpages{2526--2534}.
\bdoi{10.1016/j.jbiomech.2005.07.032}
\end{barticle}
\endbibitem

\bibitem[\protect\citeauthoryear{Pan and
  Yang}{2010}]{panSurveyTransferLearning2010}
\begin{barticle}[author]
\bauthor{\bsnm{Pan},~\bfnm{Sinno~Jialin}\binits{S.~J.}} \AND
  \bauthor{\bsnm{Yang},~\bfnm{Qiang}\binits{Q.}}
(\byear{2010}).
\btitle{A {{Survey}} on {{Transfer Learning}}}.
\bjournal{IEEE Transactions on Knowledge and Data Engineering}
\bvolume{22}
\bpages{1345--1359}.
\bdoi{10.1109/TKDE.2009.191}
\end{barticle}
\endbibitem

\bibitem[\protect\citeauthoryear{Park
  et~al.}{2008}]{parkClassificationGeneFunctions2008}
\begin{barticle}[author]
\bauthor{\bsnm{Park},~\bfnm{Changyi}\binits{C.}},
  \bauthor{\bsnm{Koo},~\bfnm{Ja-Yong}\binits{J.-Y.}},
  \bauthor{\bsnm{Kim},~\bfnm{Sujong}\binits{S.}},
  \bauthor{\bsnm{Sohn},~\bfnm{Insuk}\binits{I.}} \AND
  \bauthor{\bsnm{Lee},~\bfnm{Jae~Won}\binits{J.~W.}}
(\byear{2008}).
\btitle{Classification of Gene Functions Using Support Vector Machine for
  Time-Course Gene Expression Data}.
\bjournal{Computational Statistics \& Data Analysis}
\bvolume{52}
\bpages{2578--2587}.
\bdoi{10.1016/j.csda.2007.09.002}
\end{barticle}
\endbibitem

\bibitem[\protect\citeauthoryear{Pomann, Staicu and
  Ghosh}{2016}]{pomannTwoSampleDistributionFree2016}
\begin{barticle}[author]
\bauthor{\bsnm{Pomann},~\bfnm{Gina-Maria}\binits{G.-M.}},
  \bauthor{\bsnm{Staicu},~\bfnm{Ana-Maria}\binits{A.-M.}} \AND
  \bauthor{\bsnm{Ghosh},~\bfnm{Sujit}\binits{S.}}
(\byear{2016}).
\btitle{A {{Two Sample Distribution-Free Test}} for {{Functional Data}} with
  {{Application}} to a {{Diffusion Tensor Imaging Study}} of {{Multiple
  Sclerosis}}}.
\bjournal{Journal of the Royal Statistical Society. Series C, Applied
  statistics}
\bvolume{65}
\bpages{395--414}.
\bdoi{10.1111/rssc.12130}
\end{barticle}
\endbibitem

\bibitem[\protect\citeauthoryear{Ramsay and
  Silverman}{2002}]{ramsayAppliedFunctionalData2002}
\begin{bbook}[author]
\beditor{\bsnm{Ramsay},~\bfnm{James~O.}\binits{J.~O.}} \AND
  \beditor{\bsnm{Silverman},~\bfnm{Bernard~W.}\binits{B.~W.}}, eds.
(\byear{2002}).
\btitle{Applied {{Functional Data Analysis}}: {{Methods}} and {{Case
  Studies}}}.
\bseries{Springer {{Series}} in {{Statistics}}}.
\bpublisher{{Springer}}, \baddress{{New York, NY}}.
\bdoi{10.1007/b98886}
\end{bbook}
\endbibitem

\bibitem[\protect\citeauthoryear{Reeve, Cannings and
  Samworth}{2021}]{reeveAdaptiveTransferLearning2021}
\begin{barticle}[author]
\bauthor{\bsnm{Reeve},~\bfnm{Henry W.~J.}\binits{H.~W.~J.}},
  \bauthor{\bsnm{Cannings},~\bfnm{Timothy~I.}\binits{T.~I.}} \AND
  \bauthor{\bsnm{Samworth},~\bfnm{Richard~J.}\binits{R.~J.}}
(\byear{2021}).
\btitle{Adaptive Transfer Learning}.
\bjournal{The Annals of Statistics}
\bvolume{49}
\bpages{3618--3649}.
\bdoi{10.1214/21-AOS2102}
\end{barticle}
\endbibitem

\bibitem[\protect\citeauthoryear{Rice and
  Silverman}{1991}]{riceEstimatingMeanCovariance1991}
\begin{barticle}[author]
\bauthor{\bsnm{Rice},~\bfnm{John~A.}\binits{J.~A.}} \AND
  \bauthor{\bsnm{Silverman},~\bfnm{B.~W.}\binits{B.~W.}}
(\byear{1991}).
\btitle{Estimating the {{Mean}} and {{Covariance Structure Nonparametrically
  When}} the {{Data}} Are {{Curves}}}.
\bjournal{Journal of the Royal Statistical Society: Series B (Methodological)}
\bvolume{53}
\bpages{233--243}.
\bdoi{10.1111/j.2517-6161.1991.tb01821.x}
\end{barticle}
\endbibitem

\bibitem[\protect\citeauthoryear{Song
  et~al.}{2008}]{songOptimalClassificationTimecourse2008}
\begin{barticle}[author]
\bauthor{\bsnm{Song},~\bfnm{Joon~Jin}\binits{J.~J.}},
  \bauthor{\bsnm{Deng},~\bfnm{Weiguo}\binits{W.}},
  \bauthor{\bsnm{Lee},~\bfnm{Ho-Jin}\binits{H.-J.}} \AND
  \bauthor{\bsnm{Kwon},~\bfnm{Deukwoo}\binits{D.}}
(\byear{2008}).
\btitle{Optimal Classification for Time-Course Gene Expression Data Using
  Functional Data Analysis}.
\bjournal{Computational Biology and Chemistry}
\bvolume{32}
\bpages{426--432}.
\bdoi{10.1016/j.compbiolchem.2008.07.007}
\end{barticle}
\endbibitem

\bibitem[\protect\citeauthoryear{Staicu
  et~al.}{2012}]{staicuModelingFunctionalData2012}
\begin{barticle}[author]
\bauthor{\bsnm{Staicu},~\bfnm{Ana-Maria}\binits{A.-M.}},
  \bauthor{\bsnm{Crainiceanu},~\bfnm{Ciprian~M.}\binits{C.~M.}},
  \bauthor{\bsnm{Reich},~\bfnm{Daniel~S.}\binits{D.~S.}} \AND
  \bauthor{\bsnm{Ruppert},~\bfnm{David}\binits{D.}}
(\byear{2012}).
\btitle{Modeling {{Functional Data With Spatially Heterogeneous Shape
  Characteristics}}}.
\bjournal{Biometrics}
\bvolume{68}
\bpages{331--343}.
\bdoi{10.1111/j.1541-0420.2011.01669.x}
\end{barticle}
\endbibitem

\bibitem[\protect\citeauthoryear{Tian and
  Feng}{2022}]{tianTransferLearningHighDimensional2022}
\begin{barticle}[author]
\bauthor{\bsnm{Tian},~\bfnm{Ye}\binits{Y.}} \AND
  \bauthor{\bsnm{Feng},~\bfnm{Yang}\binits{Y.}}
(\byear{2022}).
\btitle{Transfer {{Learning Under High-Dimensional Generalized Linear
  Models}}}.
\bjournal{Journal of the American Statistical Association}
\bvolume{0}
\bpages{1--14}.
\bdoi{10.1080/01621459.2022.2071278}
\end{barticle}
\endbibitem

\bibitem[\protect\citeauthoryear{Tzeng
  et~al.}{2017}]{tzengAdversarialDiscriminativeDomain2017}
\begin{binproceedings}[author]
\bauthor{\bsnm{Tzeng},~\bfnm{Eric}\binits{E.}},
  \bauthor{\bsnm{Hoffman},~\bfnm{Judy}\binits{J.}},
  \bauthor{\bsnm{Saenko},~\bfnm{Kate}\binits{K.}} \AND
  \bauthor{\bsnm{Darrell},~\bfnm{Trevor}\binits{T.}}
(\byear{2017}).
\btitle{Adversarial {{Discriminative Domain Adaptation}}}.
In \bbooktitle{2017 {{IEEE Conference}} on {{Computer Vision}} and {{Pattern
  Recognition}} ({{CVPR}})}
\bpages{2962--2971}.
\bdoi{10.1109/CVPR.2017.316}
\end{binproceedings}
\endbibitem

\bibitem[\protect\citeauthoryear{Wainwright}{2019}]{wainwrightHighdimensionalStatisticsNonasymptotic2019}
\begin{bbook}[author]
\bauthor{\bsnm{Wainwright},~\bfnm{Martin~J.}\binits{M.~J.}}
(\byear{2019}).
\btitle{High-{{Dimensional Statistics}}: {{A Non-Asymptotic Viewpoint}}}.
\bseries{Cambridge {{Series}} in {{Statistical}} and {{Probabilistic
  Mathematics}}}.
\bpublisher{{Cambridge University Press}}, \baddress{{Cambridge}}.
\bdoi{10.1017/9781108627771}
\end{bbook}
\endbibitem

\bibitem[\protect\citeauthoryear{Wang, Chiou and
  M{\"u}ller}{2016}]{wangFunctionalDataAnalysis2016}
\begin{barticle}[author]
\bauthor{\bsnm{Wang},~\bfnm{Jane-Ling}\binits{J.-L.}},
  \bauthor{\bsnm{Chiou},~\bfnm{Jeng-Min}\binits{J.-M.}} \AND
  \bauthor{\bsnm{M{\"u}ller},~\bfnm{Hans-Georg}\binits{H.-G.}}
(\byear{2016}).
\btitle{Functional Data Analysis}.
\bjournal{Annual Review of Statistics and Its Application}
\bvolume{3}
\bpages{257--295}.
\bdoi{10.1146/annurev-statistics-041715-033624}
\end{barticle}
\endbibitem

\bibitem[\protect\citeauthoryear{Weaver, Xiao and
  Lu}{2021}]{weaverFunctionalDataAnalysis2021}
\begin{barticle}[author]
\bauthor{\bsnm{Weaver},~\bfnm{Caleb}\binits{C.}},
  \bauthor{\bsnm{Xiao},~\bfnm{Luo}\binits{L.}} \AND
  \bauthor{\bsnm{Lu},~\bfnm{Wenbin}\binits{W.}}
(\byear{2021}).
\btitle{Functional Data Analysis for Longitudinal Data with Informative
  Observation Times}.
\bjournal{Biometrics}.
\bdoi{10.1111/biom.13646}
\end{barticle}
\endbibitem

\bibitem[\protect\citeauthoryear{Weiss, Khoshgoftaar and
  Wang}{2016}]{weissSurveyTransferLearning2016}
\begin{barticle}[author]
\bauthor{\bsnm{Weiss},~\bfnm{Karl}\binits{K.}},
  \bauthor{\bsnm{Khoshgoftaar},~\bfnm{Taghi~M.}\binits{T.~M.}} \AND
  \bauthor{\bsnm{Wang},~\bfnm{DingDing}\binits{D.}}
(\byear{2016}).
\btitle{A Survey of Transfer Learning}.
\bjournal{Journal of Big Data}
\bvolume{3}
\bpages{9}.
\bdoi{10.1186/s40537-016-0043-6}
\end{barticle}
\endbibitem

\bibitem[\protect\citeauthoryear{Wu and
  M{\"u}ller}{2010}]{wuFunctionalEmbeddingClassification2010}
\begin{barticle}[author]
\bauthor{\bsnm{Wu},~\bfnm{Ping-Shi}\binits{P.-S.}} \AND
  \bauthor{\bsnm{M{\"u}ller},~\bfnm{Hans-Georg}\binits{H.-G.}}
(\byear{2010}).
\btitle{Functional Embedding for the Classification of Gene Expression
  Profiles}.
\bjournal{Bioinformatics}
\bvolume{26}
\bpages{509--517}.
\bdoi{10.1093/bioinformatics/btp711}
\end{barticle}
\endbibitem

\bibitem[\protect\citeauthoryear{Yuan and
  Cai}{2010}]{yuanReproducingKernelHilbert2010}
\begin{barticle}[author]
\bauthor{\bsnm{Yuan},~\bfnm{Ming}\binits{M.}} \AND
  \bauthor{\bsnm{Cai},~\bfnm{T.~Tony}\binits{T.~T.}}
(\byear{2010}).
\btitle{A Reproducing Kernel {{Hilbert}} Space Approach to Functional Linear
  Regression}.
\bjournal{The Annals of Statistics}
\bvolume{38}
\bpages{3412--3444}.
\bdoi{10.1214/09-AOS772}
\end{barticle}
\endbibitem

\bibitem[\protect\citeauthoryear{Zhou, Lin and
  Liang}{2018}]{zhouEfficientEstimationNonparametric2018}
\begin{barticle}[author]
\bauthor{\bsnm{Zhou},~\bfnm{Ling}\binits{L.}},
  \bauthor{\bsnm{Lin},~\bfnm{Huazhen}\binits{H.}} \AND
  \bauthor{\bsnm{Liang},~\bfnm{Hua}\binits{H.}}
(\byear{2018}).
\btitle{Efficient {{Estimation}} of the {{Nonparametric Mean}} and {{Covariance
  Functions}} for {{Longitudinal}} and {{Sparse Functional Data}}}.
\bjournal{Journal of the American Statistical Association}
\bvolume{113}
\bpages{1550--1564}.
\bdoi{10.1080/01621459.2017.1356317}
\end{barticle}
\endbibitem

\end{thebibliography}
